\documentclass{article}

\usepackage{graphicx} 
\usepackage{amsfonts}
\usepackage{amssymb} 
\usepackage{latexsym}
\usepackage{pifont}
\usepackage{amssymb,amsfonts}
\usepackage{amsthm}
\usepackage{amsmath}
\usepackage{graphics,graphicx,psfrag}
\usepackage{amscd}
\usepackage{epsfig}
\usepackage[all]{xy}
\usepackage{float}

\newtheorem{theorem}{Theorem}[section]
\newtheorem{proposition}{Proposition}[section]
\newtheorem{definition}{Definition}[section]
\newtheorem{remark}{Remark}[section]
\newtheorem{lemma}{Lemma}[section]
\newtheorem{example}{Example}[section]

\newtheorem{corollary}{Corollary}[section]

\usepackage[hang,flushmargin]{footmisc} 
\usepackage{framed}
\usepackage[usenames,dvipsnames]{color}
\usepackage{tikz}
\usepackage{xcolor}
\usepackage{parskip}
\usetikzlibrary{matrix,arrows}
\usetikzlibrary{decorations.pathmorphing}
\usepackage{enumitem}
\usepackage{setspace}
\setlist{nolistsep}

\usepackage{setspace}
\newcommand{\Modl}{ \mbox{}_R{\rm{\bf Mod}} }

\newcommand{\Cadr}{ {\rm {\bf Ch}}({\rm {\bf Mod}}_R) }
\newcommand{\Cadl}{ {\rm {\bf Ch}}(\mbox{}_R{\rm {\bf Mod}}) }
\newcommand{\Complexes}{ {\rm {\bf Ch}}(\mathcal{C}) }   

\begin{document}

\begin{center}
{\bf HOMOLOGICAL DIMENSIONS AND ABELIAN MODEL STRUCTURES ON CHAIN COMPLEXES} 
\end{center} 

\begin{center}
{\sc Marco A. P\'erez B.} \\
Universit\'e du Qu\'ebec \`a Montr\'eal. \\
D\'epartement de Math\'ematiques. \\
marco.perez@cirget.ca
\end{center}

\begin{center} April, 2014.
\end{center}


\begin{abstract} \noindent We construct Abelian model structures on the category of chain complexes over a ring $R$, from the notion homological dimensions of modules. Given an integer $n > 0$, we prove that the left modules over a ringoid $\mathfrak{R}$ with projective dimension at most $n$ form the left half of a complete cotorsion pair. Using this result we prove that there is a unique Abelian model structure on the category of chain complexes where the exact complexes are the trivial objects and the complexes with projective dimension at most $n$ form the class of trivially cofibrant objects. In \cite{Rada}, the authors construct an Abelian model structure on chain complexes, where the trivial objects are the exact complexes and the class of cofibrant objects is given by the complexes whose terms are all projective. We extend this result by finding a new Abelian model structure with the same trivial objects and where the cofibrant objects are given by the class of complexes whose terms are modules with projective dimension at most $n$. We also prove similar results concerning the flat dimension.
\end{abstract}


\section{Introduction}

In 2002, M. Hovey stablished in \cite[Theorem 2.2]{Hovey2} a correspondence to construct an Abelian model structure from two compatible and complete cotorsions pairs (See {\cite[Theorem 2.2]{Hovey2}}). We apply this correspondence in order to construct new model structures on the category of chain complexes, from certain special classes of chain complexes involving homological dimensions. In his book \cite[Section 2.3]{Hovey}, Hovey constructs a model category structure on the category $\Cadl$ of chain complexes of left $R$-modules, where the weak equivalences are the quasi-isomorphisms and the trivial cofibrations are the monomorphisms whose cokernels are projective. We call this structure the \underline{projective model structure}. We extend this model structure to any projective dimension. Specifically, we shall find a model structure with the same weak equivalences and where the trivial cofibrations are the monomorphisms whose cokernels are complexes with projective dimension $\leq n$. The method we use to construct this structure involves the category of left modules over a ringoid $\mathfrak{R}$, denoted ${\rm {\bf Mod}}(\mathfrak{R})$. If we consider the class $\mathcal{P}_n({\rm {\bf Mod}}(\mathfrak{R}))$ of $n$-projective modules over $\mathfrak{R}$ (i.e. modules with projective dimension $\leq n$), then we show that every module in $\mathcal{P}_n({\rm {\bf Mod}}(\mathfrak{R}))$ is filtered by the set $(\mathcal{P}_n({\rm {\bf Mod}}(\mathfrak{R})))^{\leq \kappa}$ of $n$-projective modules with cardinality $\leq \kappa$, where $\kappa$ is a fixed infinite regular cardinal satisfying a certain condition for $\mathfrak{R}$. This implies that $\mathcal{P}_n({\rm {\bf Mod}}(\mathfrak{R}))$ is the left half of a complete cotorsion pair. These facts in ${\rm {\bf Mod}}(\mathfrak{R})$ are generalizations of \cite[Proposition 4.1 and Theorem 4.2]{Aldrich}. Since the category $\Cadl$ of chain complexes of left $R$-modules is a particular example of ${\rm {\bf Mod}}(\mathfrak{R})$, we have that the class of $n$-projective complexes is the left half of a complete cotorsion pair, which turns out to be compatible with another complete cotorsion pair whose left half is given by the class of differential graded $n$-projective complexes. 

In \cite{Rada}, D. Bravo, E. Enochs, A. Iacob, O. Jenda and J. Rada construct an Abelian model structure on $\Cadl$, where the trivial objects are the exact complexes and the class of cofibrant objects is given by the complexes with projective terms. Among the arguments given by the authors, it is the famous theorem by I. Kaplansky on projective modules, which states that every projective module is a direct sum of countably generated projective modules. We present an extension of Kaplansky's theorem to any projective dimension, in order to prove that every degreewise (dw) $n$-projective complex (we mean a complex whose terms have projective dimension $\leq n$) is filtered by the set of chain complexes $X$ such that each module $X_k$ has a projective resolution of length $n$, where each projective term in the resolution is written as a direct sum of countably generated projective modules over a countable set. We introduce the notion of nice sub-resolutions of $n$-projective modules to show that every exact dw-$n$-projective complex is filtered by the set of exact $dw$-$n$-projective complexes with cardinality $\leq \kappa$, with $\kappa$ as above. From these results we obtain two compatible and complete cotorsion pairs from which we have a new Abelian model structure on $\Cadl$ where the quasi-isomorphisms are the weak equivalences and the monomorphisms with (exact) dw-$n$-projective cokernels are the (trivial) cofibrations. 

With respect to the flat dimension, S. Aldrich, E. Enochs, J. R. Garc\'ia Rozas and L. Oyonarte proved in \cite[Propositions 3.1 and 4.1]{Aldrich1} that every flat complex is filtered by the set of  flat complexes with cardinality $\leq \kappa$. They also proved a similar result for the class of complexes which have all their terms flat. Their arguments can be adapted to prove  analogous results for the classes of $n$-flat complexes (i.e. complexes with flat dimension $\leq n$) and degreewise $n$-flat complexes (i.e. complexes which have all their terms with flat dimension $\leq n$).


\section{Preliminaries}

Given an Abelian category $\mathcal{C}$, two classes $\mathcal{A}$ and $\mathcal{B}$ of objects in $\mathcal{C}$ form a \underline{cotorsion pair} $(\mathcal{A,B})$ if they are orthogonal to each other with respect to the first extension bi-functor ${\rm Ext}^1_{\mathcal{C}}(-,-)$, i.e. the following equalities are satisfied:
\begin{itemize}
\item[$\bullet$] $\mathcal{A} = \mbox{}^{\perp}\mathcal{B} := \{ X \in {\rm Ob}(\mathcal{C}) \mbox{ / } {\rm Ext}^1_{\mathcal{C}}(X, B) = 0 \mbox{ for every }B \in \mathcal{B} \}$.

\item[$\bullet$] $\mathcal{B} = \mathcal{A}^{\perp} := \{ X \in {\rm Ob}(\mathcal{C}) \mbox{ / } {\rm Ext}^1_{\mathcal{C}}(A, X) = 0 \mbox{ for every }A \in \mathcal{A} \}$.
\end{itemize} 
We are interested in a special class of cotorsion pairs $(\mathcal{A,B})$ from which it is possible to obtain pre-covers and pre-envelopes by the left and right halves of $(\mathcal{A,B})$, respectively. A cotorsion pair $(\mathcal{A, B})$ in $\mathcal{C}$ is \underline{complete} if for every object $X$ of $\mathcal{C}$ there exist short exact sequences $0 \longrightarrow B \longrightarrow A \longrightarrow X \longrightarrow 0$ and $0 \longrightarrow X \longrightarrow B' \longrightarrow A' \longrightarrow 0$ where $A, A' \in \mathcal{A}$ and $B, B' \in \mathcal{B}$. The map $A \longrightarrow X$ is called \underline{special $\mathcal{A}$-pre-cover} and $X \longrightarrow B'$ \underline{special $\mathcal{B}$-pre-envelope}. 

Two cotorsion pairs $(\mathcal{A}, \mathcal{B}')$ and $(\mathcal{A}', \mathcal{B})$ are \underline{compatible} if there exists a class of objects $\mathcal{W}$ in $\mathcal{C}$  such that $\mathcal{A}' = \mathcal{A} \cap \mathcal{W}$ and $\mathcal{B}' = \mathcal{B} \cap \mathcal{W}$. 

A class of objects $\mathcal{A}$ in $\mathcal{C}$ is said to be \underline{thick} if it is closed under retracts, and if for every short exact sequence $0 \longrightarrow A' \longrightarrow A \longrightarrow A'' \longrightarrow 0$, if two out of three of the terms $A'$, $A$ and $A''$ are in $\mathcal{A}$, then so is the third. \\

\begin{theorem}[Hovey's correspondence. {\cite[Theorem 2.2]{Hovey2}}] Let $\mathcal{C}$ be a bi-complete Abelian category. If $(\mathcal{A}, \mathcal{B} \cap \mathcal{W})$ and $(\mathcal{A} \cap \mathcal{W}, \mathcal{B})$ are complete cotorsion pairs in $\mathcal{C}$ and the class $\mathcal{W}$ is thick, then there exists a unique Abelian model structure on $\mathcal{C}$ such that $\mathcal{A}$ is the class of cofibrant objects, $\mathcal{B}$ is the class of fibrant objects, and $\mathcal{W}$ is the class of trivial objects. 

Conversely, if $\mathcal{C}$ is equipped with an Abelian model structure and $\mathcal{A}$, $\mathcal{B}$ and $\mathcal{W}$ denote the classes of cofibrant, fibrant and trivial objects, respectively, then $(\mathcal{A} \cap \mathcal{W}, \mathcal{B})$ and $(\mathcal{A}, \mathcal{B} \cap \mathcal{W})$ are complete cotorsion pairs. \\
\end{theorem}

Roughly speaking, a \underline{model structure} on a bi-complete category $\mathcal{C}$ is given by three classes of morphisms of $\mathcal{C}$, called cofibrations, fibrations and weak equivalences, such that it is possible to do homotopy theory in $\mathcal{C}$. A trivial or acyclic cofibration is a map which is a cofibration and a weak equivalence. Trivial fibrations are defined similarly. We do not recall the axioms defining a model structure, but a good introduction to this notion along with very detailed examples is the book \cite{Hovey} by M. Hovey. 

If we are given a model structure on an Abelian category $\mathcal{C}$, an object $X$ of $\mathcal{C}$ is said to be \underline{cofibrant} (resp. \underline{fibrant}) if the map $0 \longrightarrow X$ is a cofibration (resp. if $X \longrightarrow 0$ is a fibration). If $0 \longrightarrow X$ is a weak equivalence, then $X$ is said to be \underline{trivial} or \underline{acyclic}. 

A model structure on a bi-complete Abelian category $\mathcal{C}$ is said to be \underline{Abelian} if the following two conditions are satisfied:
\begin{itemize}
\item[$\bullet$] A map $f$ is a (trivial) cofibration if, and only if, $f$ is a monomorphism and ${\rm CoKer}(f)$ is cofibrant (and trivial).

\item[$\bullet$] A map $g$ is a (trivial) fibration if, and only if, $g$ is an epimorphism and ${\rm Ker}(g)$ is fibrant (and trivial). 
\end{itemize}

The Eklof and Trlifaj Theorem is probably one of the most effective methods to obtain complete cotorsion pairs. It states that if $\mathcal{S}$ is a set of objects in a Grothendieck category then the cotorsion pair $(\mbox{}^\perp(\mathcal{S}^\perp), \mathcal{S}^\perp)$ is complete (This theorem was originally proven in \cite[Theorem 10]{Eklof} for the category of modules, but it is also valid in every Grothendieck category, as specified in \cite[Remark 3.2.2 (b)]{Gobel}.). 

If $(\mathcal{A,B})$ is a cotorsion pair such that $\mathcal{B} = \mathcal{S}^\perp$ for some set $\mathcal{S}$, then $(\mathcal{A,B})$ is said to be \underline{cogenerated by $\mathcal{S}$}. In this work, every time we want to show that a cotorsion pair $(\mathcal{A,B})$ is complete, we shall find a set $\mathcal{S}$ such that every object in $\mathcal{A}$ is $\mathcal{S}$-filtered. Under this conditions, the pair $(\mathcal{A,B})$ is cogenerated by $\mathcal{S}$. 

A \underline{transfinite composition} in a cocomplete Abelian category $\mathcal{C}$ is a morphism of the form $f : F_0 \longrightarrow {\rm CoLim}_{\alpha < \lambda}(F_\alpha)$, where $F : [\lambda] \longrightarrow \mathcal{C}$ is a colimit preserving functor and $\lambda$ is an ordinal. The morphism $f$ is also known as the transfinite composition of the morphisms $F_\alpha \longrightarrow F_{\alpha + 1}$ for every $\alpha + 1 < \lambda$. If in addition, all the morphisms $F_\alpha \longrightarrow F_{\alpha + 1}$ are monic with cokernel in some class $\mathcal{S}$, the $F_0 \stackrel{f}\longrightarrow {\rm CoLim}_{\alpha < \lambda} (F_\alpha)$ is called a \underline{transfinite extension} of $F_0$ by $\mathcal{S}$. If $F_0 \in \mathcal{S}$ as well, the colimit ${\rm CoLim}_{\alpha < \lambda} (F_\alpha)$ is called a \underline{transfinite extension of $\mathcal{S}$}. 

If $X$ is an object of $\mathcal{C}$ such that $X = {\rm CoLim}_{\alpha < \lambda} (F_\alpha)$ is a transfinite extension of $\mathcal{S}$, then $X$ is said to be \underline{$\mathcal{S}$-filtered} and the family $(F_\alpha \mbox{ : } \alpha < \lambda)$ is said to be an \underline{$\mathcal{S}$-filtration}. \\

\begin{definition} Let $\mathcal{A}$ be a class of objects in an Abelian category $\mathcal{C}$. We shall say that an object $X$ of $\mathcal{C}$ is a left $n$-$\mathcal{A}$-object if there is an exact sequence $0 \longrightarrow A_n \longrightarrow A_{n-1} \longrightarrow \cdots  \longrightarrow A_0 \longrightarrow X \longrightarrow 0$ with $A_k \in \mathcal{A}$ for every $0 \leq k \leq n$. Right $n$-$\mathcal{A}$-objects are defined dually. \\
\end{definition}

We recall that a cotorsion pair $(\mathcal{A,B})$ in $\mathcal{C}$ is \underline{hereditary} if ${\rm Ext}^i_{\mathcal{C}}(A,B) = 0$ for every $A \in \mathcal{A}$, $B \in \mathcal{B}$ and $i > 1$. The following result is easy to prove. \\

\begin{lemma}\label{cerradoporexts} Let $(\mathcal{A,B})$ be a hereditary cotorsion pair in an Abelian category $\mathcal{C}$. Then the classes of left $n$-$\mathcal{A}$-objects and right $n$-$\mathcal{B}$-objects are both closed under extensions. 
\end{lemma}

The following lemma is proven in the category of left $R$-modules in the reference given below, but the arguments appearing there carry over to every Grothendieck category. \\

\begin{lemma}[Eklof Lemma. {\cite[Lemma 3.1.2]{Gobel}}] Let $X$ and $Y$ be two objects of a cocomplete Abelian category $\mathcal{C}$, and suppose $X = {\rm CoLim}_{\alpha < \lambda} (X_\alpha)$, where $(X_\alpha : \alpha < \lambda)$ is a transfinite extension of $\mbox{}^\perp\{ Y \}$. Then $X \in \mbox{}^\perp\{ Y \}$. 
\end{lemma}

As a consequence, we have the following proposition. \\

\begin{proposition}\label{gentransex} Let $(\mathcal{A,B})$ be a cotorsion pair in an Abelian category $\mathcal{C}$. If $\mathcal{S} \subseteq \mathcal{A}$ is a set objects in $\mathcal{C}$ such that every object $A \in \mathcal{A}$ is a transfinite extension of $\mathcal{S}$, then $(\mathcal{A,B})$ is cogenerated by $\mathcal{S}$. 
\end{proposition}


\section{Cotorsion pairs from projective dimensions of left modules over a ringoid}

In this section, we shall study some facts concerning the projective dimension of objects in an Abelian category $\mathcal{C}$. We shall focus in the particular case where $\mathcal{C}$ is the category ${\rm {\bf Mod}}(\mathfrak{R})$ of modules over a ringoid $\mathfrak{R}$, for which we prove the class of objects with projective dimension bounded by some $n \geq 0$ is the left half of a cotorsion pair cogenerated by a set. This was initially proven by S. T. Aldrich and coauthors in \cite{Aldrich} in the case $\mathcal{C}$ is the category of left $R$-modules $\Modl$. We do some modifications to their arguments to generalize their result to ${\rm {\bf Mod}}({\mathfrak{R}})$. The advantage in doing this lies in the fact that both $\Modl$ and the category $\Cadl$ of complexes of left $R$-moules are particular examples of categories of modules over a ringoid. 

We recall some notations in the category $\Complexes$ of complexes over an Abelian category $\mathcal{C}$. Given a chain complex $X = (X_m)_{m \in \mathbb{Z}}$ with boundary or differential maps $\partial^X_m : X_m \longrightarrow X_{m-1}$, the object $Z_m(X) := {\rm Ker}(\partial^X_m)$ is called the \underline{$m$-cycle} of $X$, and $B_m (X) := {\rm Im}(\partial^X_{m+1})$ the \underline{$m$-boundary} of $X$. \\

\begin{example} Let $C \in {\rm Ob}(\mathcal{C})$. The \underline{$m$th sphere complex centred at $C$}, denoted $S^m(C)$, is defined by $(S^m(C))_m := C$ and $(S^m(C))_k := 0$ if $k \neq m$. The \underline{$m$th disk complex centred at $C$}, denoted $D^m(C)$, is defined by $(D^m(C))_k := C$ if $k = m$ or $m-1$, and $0$ otherwise, where the boundary map $\partial^{D^m(C)}_m$ is the identity ${\rm id}_C$. 
\end{example}

From now on, left $n$-$\mathcal{P}_0(\mathcal{C})$-objects and right $n$-$\mathcal{I}_0(\mathcal{C})$-objects shall be called \underline{$n$-projective} and \underline{$n$-injective objects}, respectively. Let $X$ be an object in an Abelian category $\mathcal{C}$. Let ${\rm pd}(X)$ and ${\rm id}(X)$ denote the projective and injective dimension of $X$, respectively. Note that $X$ is $n$-projective (resp. $n$-injective) if, and only if, ${\rm pd}(X) \leq n$ (resp. ${\rm id}(X) \leq n$). We shall denote by $\mathcal{P}_n(\mathcal{C})$ and $\mathcal{I}_n(\mathcal{C})$ the classes of $n$-projective and $n$-injective objects of $\mathcal{C}$, respectively. \\

\begin{example}[{\cite[Theorem 4.2]{Aldrich}}] \label{nprojpair} Let $R$ be an associative ring with unity and $\kappa$ be an infinite cardinal such that $\kappa \geq {\rm Card}(R)$. The class $\mathcal{P}_n := \mathcal{P}_n(\Modl)$ of $n$-projective modules is the left half of a cotorsion pair $\left(\mathcal{P}_n, (\mathcal{P}_n)\mbox{}^\perp\right)$ cogenerated by the set $(\mathcal{P}_n)^{\leq \kappa} := \{ M \in \mathcal{P}_n \mbox{ {\rm :} } {\rm Card}(M) \leq \kappa \}$. 
\end{example}

Using induction and the fact that the class of exact chain complexes is thick, the following result follows. \\

\begin{lemma} \label{lema1} Let $\mathcal{C}$ be an Abelian category. 
\begin{itemize}
\item[{\bf (1)}] If $0 \longrightarrow A_n \longrightarrow \cdots \longrightarrow A_0 \stackrel{f_0}\longrightarrow X \longrightarrow 0$ is an exact sequence in $\Complexes$ such that $A_i$ is exact for every $0 \leq i \leq n$, then so is $X$. 

\item[{\bf (2)}] If $0 \longrightarrow Y \stackrel{g^0}\longrightarrow B^0 \longrightarrow \cdots \longrightarrow B^n \longrightarrow 0$ is an exact sequence in $\Complexes$ such that $B^i$ is exact for every $0 \leq i \leq n$, then so is $Y$. 
\end{itemize} 
\end{lemma}

\begin{lemma} Consider a short exact sequence $S = 0 \longrightarrow Y \stackrel{f}\longrightarrow Z \stackrel{g}\longrightarrow X \longrightarrow 0$ in $\Complexes$. 
\begin{itemize}
\item[{\bf (1)}] The sequence $0 \longrightarrow Z_m(Y) \longrightarrow Z_m(Z) \longrightarrow Z_m(X) \longrightarrow 0$ is exact if $Y$ is an exact complex.

\item[{\bf (2)}] The sequence $0 \longrightarrow \frac{Y_m}{B_m(Y)} \longrightarrow \frac{Z_m}{B_m(Z)} \longrightarrow \frac{X_m}{B_m(X)} \longrightarrow 0$ is exact if $X$ is an exact complex. 
\end{itemize}
\end{lemma}
\begin{proof} We only show the first statement for $\Modl$, by \cite[Theorem 3, page 204]{MacLane}. Let $Z_m(f) : Z_m(Y) \longrightarrow Z_m(Z)$ be the homomorphism induced by the universal property of kernels, which is given by $y \mapsto f_m(y)$ for every $y \in Z_m(Y)$. The homomorphism $Z_m(g) : Z_m(Y) \longrightarrow Z_m(X)$ is defined similarly. It is easy to check that $Z_m(f)$ is monic and that ${\rm Ker}(Z_m(g)) = {\rm Im}(Z_m(f))$. These facts do not depend on the exactness of $Y$. Let $x \in Z_m(X)$. There exists $z \in Z_m$ such that $x = g_m(z)$. We have $g_{m-1} \circ \partial^Z_m(z) = \partial^X_m \circ g_m(z) = 0$. Since the sequence $0 \longrightarrow Y_{m-1} \longrightarrow Z_{m-1} \longrightarrow X_{m-1} \longrightarrow 0$ is exact, there exists $y \in Y_{m-1}$ such that $\partial^Z_{m}(z) = f_{m-1}(y)$. Then $f_{m-2} \circ \partial^Y_{m-1}(y) = \partial^Z_{m-1} \circ f_{m-1}(y) = 0$ and so $\partial^Y_{m-1}(y) = 0$ since $f_{m-2}$. By the exactness of $Y$, there exists $y' \in Y_m$ such that $y = \partial^Y_{m}(y')$. Hence $\partial^Z_m(z - f_m(y')) = 0$ and $g_m(z - f_m(y')) = x$. 
\end{proof}

Using the previous lemma along with the induction principle, we obtain the following result. \\

\begin{lemma} \label{exactofexact} Let $0 \longrightarrow A_n \stackrel{f_n}\longrightarrow A_{n-1} \longrightarrow \cdots \longrightarrow A_1 \stackrel{f_1}\longrightarrow A_0 \longrightarrow 0$ be an exact sequence in $\Complexes$ of exact chain complexes. Then for every $m \in \mathbb{Z}$, the $m$th cycles $Z_m(A_i)$ form the following exact sequence in $\mathcal{C}$: \[ 0 \longrightarrow Z_m (A_n) \longrightarrow Z_m (A_{n-1}) \longrightarrow \cdots \longrightarrow Z_m (A_1) \longrightarrow Z_m (A_0) \longrightarrow 0. \]
\end{lemma}

Statement {\bf (2)} of the following proposition is proven in \cite[Theorem 3.1.3]{Garcia}. We give here a different argument.  \\

\begin{proposition} \label{npicomplexes} Let $\mathcal{C}$ be an Abelian category and $n$ be a positive integer. 
\begin{itemize}
\item[{\bf (1)}] Assume $\mathcal{C}$ has enough projective objects. A chain complex $X$ is $n$-projective if, and only if, $X$ is exact and each $Z_m(X)$ is an $n$-projective object of $\mathcal{C}$. 

\item[{\bf (2)}] Assume $\mathcal{C}$ has enough injective objects. A chain complex $Y$ is $n$-injective if, and only if, $Y$ is exact and each $Z_m(Y)$ is an $n$-injective object of $\mathcal{C}$.
\end{itemize} 
\end{proposition}
\begin{proof} Let $X$ be an exact complex with $n$-projective cycles. Consider a partial projective resolution $0 \longrightarrow K \longrightarrow P_{n-1} \longrightarrow \cdots \longrightarrow P_0 \longrightarrow X \longrightarrow 0$. Note $K$ is exact by Lemma \ref{lema1}. Notice also that $Z_m(P_i)$ is projective for every $0 \leq i \leq n-1$ and every $m \in \mathbb{Z}$. It follows by Lemma \ref{exactofexact} that $Z_m(K) \in \Omega^n(Z_m(X))$ (Where $\Omega^n(Z_m(X))$ denotes the class of $n$th syzygies of $Z_m(X)$.). Hence, $Z_m(K) \in \mathcal{P}_0(\mathcal{C})$ since $Z_m(X)$ is $n$-projective. The converse follows similarly. 
\end{proof}

We shall see that the class $\mathcal{P}_n(\Cadl)$ of $n$-projective chain complexes is the left half of a complete cotorsion pair. This, along with Example \ref{nprojpair}, is a consequence of a more general result in the category ${\rm {\bf Mod}}(\mathfrak{R})$ of left modules over a ringoid $\mathfrak{R}$. We begin by studying the projective dimension of objects in this category. 

A small pre-additive category $\mathfrak{R}$ is called \underline{ringoid} or \underline{ring with many objects}. Then we have a composition law  
\begin{align*}
{\rm Hom}_{\mathfrak{R}}(b,c) \otimes {\rm Hom}_{\mathfrak{R}}(a,b) & \longrightarrow {\rm Hom}_{\mathfrak{R}}(a,c) \\
y \otimes x & \mapsto y \circ x
\end{align*} 
for every $a, b, c \in {\rm Ob}(\mathfrak{R})$, and a unit ${\rm id}_a \in  {\rm Hom}_{\mathfrak{R}}(a,a)$ for every $a \in {\rm Ob}(\mathfrak{R})$. The composition law defines a ring structure on ${\rm Hom}_{\mathfrak{R}}(a,a)$ for every $a \in {\rm Ob}(\mathfrak{R}) $. Moreover, the Abelian group ${\rm Hom}_{\mathfrak{R}}(a,b)$ has the structure of a bimodule, with a left action by ${\rm Hom}_{\mathfrak{R}}(b,b)$ and a right action by ${\rm Hom}_{\mathfrak{R}}(a,a)$. \\

\begin{example}\label{examplesringoids} \
\begin{itemize}
\item[{\bf (1)}] Every ring $R$ can be regarded as a ringoid $\mathfrak{R}$ having a single object $\star$ if we put ${\rm Hom}_{\mathfrak{R}}(\star,\star) = R$. 

\item[{\bf (2)}] We shall denote by $\mathfrak{S}$ the ringoid generated by the following infinite graph
\[ \begin{tikzpicture}
\matrix (m) [matrix of math nodes, row sep=1.5em, column sep=3em]
{ \cdots & 2 & 1 & 0 & -1 & -2 & \cdots \\ };
\path[->]
(m-1-1) edge (m-1-2) (m-1-2) edge node[above] {$\partial_2$} (m-1-3) (m-1-3) edge node[above] {$\partial_1$} (m-1-4) (m-1-4) edge node[above] {$\partial_0$} (m-1-5) (m-1-5) edge node[above] {$\partial_{-1}$} (m-1-6) (m-1-6) edge (m-1-7);
\path[->,font=\scriptsize]
(m-1-2) edge[loop above] node[auto] {$e_2$} (m-1-2)
(m-1-3) edge[loop above] node[auto] {$e_3$} (m-1-3)
(m-1-4) edge[loop above] node[auto] {$e_4$} (m-1-4)
(m-1-5) edge[loop above] node[auto] {$e_5$} (m-1-5)
(m-1-6) edge[loop above] node[auto] {$e_6$} (m-1-6);
\end{tikzpicture} \]
together with the relation $\partial_n \circ \partial_{n+1} = 0$ for $n \in \mathbb{Z}$. We have ${\rm Ob}(\mathfrak{S}) = \mathbb{Z}$ and \[ {\rm Hom}_{\mathfrak{S}}(i, j) = \left\{ \begin{array}{ll}  \left< e_i \right> := \{ m \cdot e_i \mbox{ {\rm : }} m \in \mathbb{Z} \} &  \mbox{ if } j = i, \\ \left<\partial_i\right> := \{ m \cdot \partial_i \mbox{ {\rm : }} m \in \mathbb{Z} \} & \mbox{ if } j = i + 1, \\ 0 & \mbox{ otherwise}. \end{array} \right. \] 
\end{itemize} 
\end{example}

\

A \underline{(left) module $M$ over a ringoid $\mathfrak{R}$} is an additive functor $M : \mathfrak{R} \longrightarrow {\rm {\bf Ab}}$, where ${\rm {\bf Ab}}$ is the category of Abelian groups. A \underline{map of (left) $\mathfrak{R}$-modules} is a natural transformation $M \longrightarrow N$. \\

\begin{example} Given a ringoid $\mathfrak{R}$ and $a \in {\rm Ob}(\mathfrak{R}) $, the covariant functor ${\rm Hom}_{\mathfrak{R}}(a, -) : \mathfrak{R} \rightarrow {\rm {\bf Ab}}$ is a (left) $\mathfrak{R}$-module. \\
\end{example}

The category ${\rm {\bf Mod}}(\mathfrak{R})$ of left modules over $\mathfrak{R}$ is defined as the category $[\mathfrak{R}, {\rm {\bf Ab}}]$ of additive functors $\mathfrak{R} \longrightarrow {\rm {\bf Ab}}$. Note that ${\rm {\bf Mod}}(\mathfrak{R})$ is Abelian and bi-complete since ${\rm {\bf Ab}}$ is. \newpage

\begin{example}[Particular examples of ${\rm {\bf Mod}}(\mathfrak{R})$] \label{particulares} \
\begin{itemize}
\item[{\bf (1)}] If $\mathfrak{R}$ is the ringoid of Example \ref{examplesringoids} {\bf (1)}, then ${\rm {\bf Mod}}(\mathfrak{R})$ is the category $\Modl$ of left $R$-modules. 

\item[{\bf (2)}] Recall that the \underline{tensor product} $\mathcal{C} \otimes \mathcal{D}$ of two pre-additive categories $\mathcal{C}$ and $\mathcal{D}$ is the pre-additive category defined by putting ${\rm Ob}(\mathcal{C}\otimes \mathcal{D}) ={\rm Ob}(\mathcal{C}) \times {\rm Ob}(\mathcal{D})$ and ${\rm Hom}_{\mathcal{C} \otimes \mathcal{D}}((C, D), (C', D')) = {\rm Hom}_{\mathcal{C}}(C, C') \otimes {\rm Hom}_{\mathcal{D}}(D, D')$. If $\mathcal{C}$, $\mathcal{D}$ and $\mathcal{E}$ are pre-additive categories then we have a canonical isomorphism of pre-additive categories $[\mathcal{C}\otimes \mathcal{D}, \mathcal{E}]\simeq [\mathcal{C}, [\mathcal{D}, \mathcal{E}]]$. In particular, if $\mathcal{K}$ and $\mathfrak{R}$ are ringoids, then $[\mathcal{K}, \mathbf{Mod}(\mathfrak{R})] =  [\mathcal{K}, [\mathfrak{R}, \mathbf{Ab}]] \simeq [\mathcal{K}\otimes \mathfrak{R}, \mathbf{Ab}] =\mathbf{Mod}(\mathcal{K}\otimes \mathfrak{R})$. If we consider a ring $R$ and the ringoid $\mathfrak{S}$ defined in Example \ref{examplesringoids} {\bf (2)}, then we have an isomorphism $[\mathfrak{S}, \mathbf{Mod}(R)] \simeq \mathbf{Mod}(\mathfrak{S}\otimes R)$ of additive categories. This means that a chain complex of $R$-modules is a module over the ringoid $\mathfrak{S}\otimes R$, i.e. $\Cadl \cong {\rm {\bf Mod}}(\mathfrak{S} \otimes R)$. \\
\end{itemize}
\end{example}

\begin{remark} A sequence $\cdots \longrightarrow M_1 \longrightarrow M_0 \longrightarrow M_{-1} \longrightarrow \cdots$ of left $\mathfrak{R}$-modules is exact if the sequence $\cdots \longrightarrow M_1(a) \longrightarrow M_0(a) \longrightarrow M_{-1}(a) \longrightarrow \cdots$ of Abelian groups is exact for every $a \in {\rm Ob}(\mathfrak{R})$. \\
\end{remark}

\begin{definition} We shall say that an element $x\in M(a)$ is \underline{homogenous of grade} $a$ and we shall write $a = |x|$. \\
\end{definition}

If $M$ is a left module over $\mathfrak{R}$, then the map ${\rm Hom}_{\mathfrak{R}}(a,b) \longrightarrow {\rm Hom}_{{\rm {\bf Ab}}}(M(a), M(b))$ of Abelian groups defined by $M$ induces a multiplication 
\begin{align*}
{\rm Hom}_{\mathfrak{R}}(a, b) \otimes M(a) & \longrightarrow M(b) \\
(r, x) & \mapsto r \cdot x := M(r)(x)
\end{align*} 
for every $a, b \in {\rm Ob}(\mathfrak{R})$. The product of $r \in {\rm Hom}_{\mathfrak{R}}(a, b)$ by $x \in M(a)$ is an element $r \cdot x \in M(b)$. \\

\begin{definition} We shall say that a linear combination of homogenous elements $y = \sum_{i\in I} r_i \cdot x_i$ is \underline{admissible} if $y$ is homogenous and $r_i \in {\rm Hom}_{\mathfrak{R}}(|x_i|, |y|)$ for every $i\in I$. We shall accept infinite combinations in the case where $r_i = 0$ for all but finitely many $i \in I$. \\
\end{definition}

\begin{definition} If $M$ is a left $\mathfrak{R}$-module we shall say that a family $N = \{ N(a) \mbox{ : } a \in {\rm Ob}(\mathfrak{R}) \}$ of subgroups $N(a) \subseteq M(a)$ is a \underline{sub-module} of $M$ if $x\in N(a)$ implies $r \cdot x \in N(b)$ for every $r \in {\rm Hom}_{\mathfrak{R}}(a,b)$. \\
\end{definition}

\begin{remark} Note that the family $\{ N(a) \mbox{ : } a \in {\rm Ob}(\mathfrak{R}) \}$ in the previous definition defines a functor $N : \mathfrak{R} \longrightarrow {\rm {\bf Ab}}$, where $N(r)$ is the restriction $M(r)|_{N(a)}$ for every map $r : a \longrightarrow b$. Conversely, if $N$ is a sub-functor of $M$, then $\{ N(a) \mbox{ : } a \in {\rm Ob}(\mathfrak{R}) \}$ is a sub-module of $M$. \\
\end{remark}

We want to construct for every $n$-projective left $\mathfrak{R}$-module, a transfinite extension of \textquotedblleft small\textquotedblright \ $n$-projective $\mathfrak{R}$-modules. The construction of these transfinite extension in the case $\mathcal{C} = \Modl$ is based on a method, probably first introduced by S. T. Aldrich and coauthors in \cite{Aldrich}, known as the {\it zig-zag procedure}. We shall explain how to adapt this procedure to the category ${\rm {\bf Mod}}(\mathfrak{R})$. 

We need to introduce some notation and recall the notion of free $\mathfrak{R}$-modules. \\

\begin{definition} If $M$ is a $\mathfrak{R}$-module we shall say that a sub-module $N \subseteq M$ is \underline{generated} by a family $\{ x_i \}_{i \in I}$ of homogenous elements if $N$ is the smallest sub-module of $M$ which contains all the elements $x_i$. \\
\end{definition}

Let $N \subseteq M$ be a sub-module generated by a family $\{ x_i \}_{i \in I}$ of homogenous elements of $M$. Then an element $x \in M(a)$ belongs to $N(a)$ if and only if it is an admissible linear combination $x = \sum_{i \in I} r_i \cdot x_i$ (where $r_i = 0$ for all but finitely many $i\in I$). \\

\begin{definition} We shall say that a family $\{ x_i \}_{i\in I}$ of homogenous elements of $M$ is a \underline{basis} of $M$ if  every homogenous element $x \in M$ can be written uniquely as an admissible linear combination $x = \sum_{i \in I} r_i \cdot x_i$. We shall say that $M$ is \underline{free} if it admits a basis. \\
\end{definition}

If $\mathfrak{R}$ is a ringoid and $M$ is a left $\mathfrak{R}$-module, then it follows from Yoneda's Lemma that for every $a \in {\rm Ob}(\mathfrak{R})$ and every $x \in M(a)$ there
is a unique map of $\mathfrak{R}$-modules $\alpha : {\rm Hom}_{\mathfrak{R}}(a, -) \longrightarrow M$ such that $x = \alpha_a ({\rm id}_a)$. More generally, if $\{ a_i \}_{i\in I}$ be a family of objects of $\mathfrak{R}$, let us put \[ \left< a_i \mbox{ : } i \in I \right> := \bigoplus_{i \in I} {\rm Hom}_{\mathfrak{R}}(a_i, -) \] and $[i] := (u_i)_{a_i}({\rm id}_{a_i})$, where $u_i: {\rm Hom}_{\mathfrak{R}}(a_i, -) \longrightarrow \left< a_i \mbox{ : } i \in I \right>$ is the inclusion. Consider a family of elements $\{ x_i \}_{i \in I}$ in $\prod_{i \in I} M(a_i)$. Then for each $i \in I$ we can write $x_i = (\alpha^i)_{a_i}({\rm id}_{a_i})$, where $\alpha^i : {\rm Hom}_{\mathfrak{R}}(a_i,-) \longrightarrow M$ is a map of left $\mathfrak{R}$-modules. Since ${\rm {\bf Mod}}(\mathfrak{R})$ is cocomplete, there exists a unique map $f : \left< a_i \mbox{ : } i \in I \right> \rightarrow M$ such that the following triangle commutes: \\
\[ \begin{tikzpicture}
\matrix (m) [matrix of math nodes, row sep=0.75cm, column sep=1.5cm]
{ {\rm Hom}_{\mathcal{R}}(a_i, -) & \left< a_i \mbox{ : } i \in I \right> \\ & M \\ };
\path[->]
(m-1-1) edge node[above] {$u_i$} (m-1-2) edge node[below,sloped] {$\alpha^i$} (m-2-2);
\path[dotted,->]
(m-1-2) edge node[right] {$\exists \mbox{! } \alpha$} (m-2-2);
\end{tikzpicture} \]
Note that $\alpha_{a_i}([i]) = \alpha_{a_i} \circ (u_i)_{a_i}({\rm id}_{a_i}) = (\alpha^i)_{a_i}({\rm id}_{a_i}) = x_i$ for every $i \in I$. It follows that the $\mathfrak{R}$-module $\left< a_i \mbox{ : } i \in I \right>$ is freely generated by the elements $[i]$ of grade $a_i$ for $i \in I$. \\

\begin{definition} The family  $\{ a_i \}_{i \in I}$ is defining  a map $| - | : I \longrightarrow {\rm Ob}(\mathfrak{R})$ if we put $| i | = a_i$ for $i \in I$. We shall say that the set $I$ equipped with the map $| - | : I \longrightarrow {\rm Ob}(\mathfrak{R})$ is \underline{$\mathfrak{R}$-graded}. \\
\end{definition}

If $I$ is an $\mathfrak{R}$-graded set, then the $\mathfrak{R}$-module $\left< I \right> = \bigoplus_{i \in I} {\rm Hom}_{\mathfrak{R}}(| i |, -)$ is freely generated by elements $[i]$ of grade $|i|$ for $i\in I$.  \\

\begin{proposition} An $\mathfrak{R}$-module $M$ is free if, and only if, it is isomorphic to a coproduct of $\mathfrak{R}$-modules ${\rm Hom}_{\mathfrak{R}}(a_i, -)$ for a family $\{ a_i \}_{i \in I}$ of objects of $\mathfrak{R}$, $M \cong \bigoplus_{i \in I} {\rm Hom}_{\mathfrak{R}}(a_i, -)$.
\end{proposition}
\begin{proof} The implication ($\Longleftarrow$) follows by the comments above. Now suppose $M$ is a free left module over $\mathfrak{R}$ admitting a basis $\{ x_i \}_{i \in I}$. Consider the natural transformation $\alpha : \left< a_i \mbox{ : } i \in I \right> \longrightarrow M$ given above, where $x_i \in M(a_i)$. We check that $\alpha$ is a natural isomorphism, i.e. $\alpha_b : \bigoplus_{i \in I} {\rm Hom}_{\mathfrak{R}}(a_i, b) \longrightarrow M(b)$ is an isomorphism for every $b \in {\rm Ob}(\mathfrak{R})$. Let $x \in M(b)$. We can write $x$ as a unique admissible linear combination $x = \sum_{i \in I} r_i \cdot x_i$, where $r_i \in {\rm Hom}_{\mathfrak{R}}(a_i, b)$ for every $i \in I$. Since $\alpha$ is a natural transformation, we have that $M(r_i) \circ \alpha_{a_i} = \alpha_b \circ \left< a_i \mbox{ : } i \in I \right>(r_i)$. Then 
\begin{align*}
x & = \sum_{i \in I} M(r_i)(x_i) = \sum_{i \in I} M(r_i) (\alpha_{a_i}([i])) = \sum_{i \in I} \alpha_b \circ \left< a_i \mbox{ : } i \in I \right>(r_i)([i]) \\
& = \sum_{i \in I} \alpha_b (r_i \cdot [i]) = \alpha_b \left( \sum_{i \in I} r_i \cdot [i] \right),  \mbox{ where $\sum_{i \in I} r_i \cdot [i]$ is unique}
\end{align*}
\end{proof}

\begin{corollary}\label{projmodfrak} Every free (left) $\mathfrak{R}$-module is projective.
\end{corollary}
\begin{proof} Since the direct sum of projective objects is projective, by the previous proposition it suffices to show that ${\rm Hom}_{\mathfrak{R}}(a,-)$ is projective in ${\rm {\bf Mod}}(\mathfrak{R})$, i.e. that the functor ${\rm Hom}_{{\rm {\bf Mod}}(\mathfrak{R})}({\rm Hom}_{\mathfrak{R}}(a,-),-) : {\rm {\bf Mod}}(\mathfrak{R}) \rightarrow {\rm {\bf Ab}}$ is exact. Suppose we are given a short exact sequence $0 \longrightarrow M' \longrightarrow M \longrightarrow M'' \longrightarrow 0$ in ${\rm {\bf Mod}}(\mathfrak{R})$. Then the sequence \\
\[ \begin{tikzpicture}
\matrix (m) [matrix of math nodes, row sep=1em, column sep=1.5em]
{ 0 & {\rm Hom}_{{\rm {\bf Mod}}(\mathfrak{R})}({\rm Hom}_{\mathfrak{R}}(a,-),M') & {\rm Hom}_{{\rm {\bf Mod}}(\mathfrak{R})}({\rm Hom}_{\mathfrak{R}}(a,-),M) \\ & {\rm Hom}_{{\rm {\bf Mod}}(\mathfrak{R})}({\rm Hom}_{\mathfrak{R}}(a,-),M'') & 0\phantom{foiasdnfoinoifaodfsdfsdfs} \\ };
\path[overlay,->, font=\scriptsize, >=latex]
(m-1-1) edge (m-1-2) (m-1-2) edge (m-1-3)
(m-1-3) edge [out=355, in=175] (m-2-2)
(m-2-2) edge (m-2-3);
\end{tikzpicture} \]
is exact since it is isomorphic to $0 \longrightarrow M'(a) \longrightarrow M(a) \longrightarrow M''(a) \longrightarrow 0$, by Yoneda's Lemma. 
\end{proof}

\

\begin{proposition}[Eilenberg's Trick for modules over a ringoid] \label{free resolution} For every projective $\mathfrak{R}$-module $P$, there exists a free $\mathfrak{R}$-module $F$ together with an isomorphism $P \oplus F \cong F$. 
\end{proposition}

As a consequence of the Eilenberg's Trick, we have: \\

\begin{corollary} Every $n$-projective (left) $\mathfrak{R}$-module has a free left resolution of length $n$. \\
\end{corollary}

\begin{definition}\label{kappasmall} Let $\kappa$ be an infinite regular cardinal strictly greater than the cardianlity of ${\rm Hom}_{\mathfrak{R}}(a, b)$ for every $a, b \in {\rm Ob}(\mathfrak{R})$. We shall say that a $\mathfrak{R}$-module $M$ is \underline{$\kappa$-small} if ${\rm Card}(M(a)) \leq \kappa$, for every $a \in {\rm Ob}(\mathfrak{R})$. \\
\end{definition}

\begin{remark}[$\kappa$-small modules and complexes] Let $\mathfrak{R}$ be a ringoid and $\kappa$ be an infinite regular cardinal such that $\kappa > {\rm Card}({\rm Hom}_{\mathfrak{R}}(a,b))$ for every $a, b \in {\rm Ob}(\mathfrak{R})$. Given a class $\mathcal{X}$ of left modules over $\mathfrak{R}$, we shall denote by $\mathcal{X}^{\leq \kappa}$ the set of $\kappa$-small left modules. 
\begin{itemize}
\item[{\bf (1)}] If $\mathfrak{R}$ is the ringoid of Example \ref{examplesringoids} {\bf (1)}, then we have that $\kappa > {\rm Card}(R)$ and that a left $R$-module $M$ is $\kappa$-small if, and only if, ${\rm Card}(M) \leq \kappa$.  

\item[{\bf (2)}] We know $\Cadl$ is equivalent to the category of left modules over the ringoid $\mathfrak{S} \otimes R$, with $\mathfrak{S}$ as in Example \ref{examplesringoids} {\bf (2)}. It is not hard to see that $\kappa > {\rm Card}(R)$ and that the following conditions are equivalent for every chain complex $X$ in $\Cadl$: 
	\begin{itemize}
	\item[{\bf (a)}] $X$ is $\kappa$-small. 
	\item[{\bf (b)}] ${\rm Card}(X_m) \leq \kappa$ for every $m \in \mathbb{Z}$.
	\item[{\bf (c)}] $\sum_{m \in \mathbb{Z}} {\rm Card}(X_m) \leq \kappa$. 
	\end{itemize}
\end{itemize}

We shall say for the rest of this work that a left $R$-module $M$ is $\kappa$-small if ${\rm Card}(M) \leq \kappa$, and that a chain complex $X$ in $\Cadl$ is $\kappa$-small if each $X_m$ is a $\kappa$-small module, where $\kappa$ is an (infinite) regular cardinal satisfying $\kappa > {\rm Card}(R)$. \\
\end{remark}

\begin{lemma}[Generalization of {\cite[Proposition 4.1]{Aldrich}}]\label{lemaoyogen} Let $\kappa$ be an infinite regular cardinal as in the previous definition. Let $M$ be a $n$-projective $\mathfrak{R}$-module. Then for every homogeneous element $x \in M(a)$ there exists a $\kappa$-small sub-module $N \hookrightarrow M$ with $x \in N(a)$ such that the $\mathfrak{R}$-modules $N$ and $M/N$ are $n$-projective.  
\end{lemma}
\begin{proof} By the previous corollary, we start with a free resolution of $M$, say \[ 0 \longrightarrow \left< I_n \right> \stackrel{\partial_n}\longrightarrow \left< I_{n-1} \right> \longrightarrow \cdots \longrightarrow \left< I_1 \right> \stackrel{\partial_1}\longrightarrow \left< I_0 \right> \stackrel{\partial_0}\longrightarrow M \longrightarrow 0, \] where $I_k$ is an $\mathfrak{R}$-graded set for every $0 \leq k \leq n$.  

The map $(\partial_0)_a : \left< I_0 \right>(a) = \bigoplus_{i \in I_0} {\rm Hom}_{\mathfrak{R}}([i],a) \longrightarrow M(a)$ is surjective, so we can find a finite number of maps $r_{i_1} : [i_1] \longrightarrow a$, $\dots$, $r_{i_k} : [i_k] \longrightarrow a$ such that $x = (\partial_0)_a(r_{i_1} + \cdots + r_{i_k})$. Then $Z_0 = \{ i_1, \dots, i_k \}$ is a finite subset of $I_0$ such that $x \in \partial_0\left( \left< Z_0 \right> \right)$ (by abuse of notation, this shall mean that $x \in (\partial_0)_a \circ \left< Z_0 \right>(a)$). Consider the natural transformation $\partial_0 |_{\left< Z_0 \right>} : \left< Z_0 \right> \longrightarrow M$ and let $y \in {\rm Ker}\left( \partial_0|_{\left< Z_0 \right>} \right)$ of degree $b$. Since $\left< I_1 \right>(b) \longrightarrow \left< I_0 \right>(b) \longrightarrow M(b)$ is exact, there exists $y' \in \left< I_1\right>(b)$ such that $y = (\partial_1)_b(y')$. We can write $y' = \sum_{i \in Z^y_1} r_i \cdot [i]$, where $Z^y_1$ is a finite subset of $I_1$. Let $Z_1 = \bigsqcup \{ Z^y_1 \mbox{ : } y \in {\rm Ker}\left( \partial_0 |_{\left< Z_0 \right>} \right) \}$. In order to estimate the number of elements of $Z^y_1$, note that for each $y \in {\rm Ker}\left( \partial_0 |_{\left< Z_0 \right>} \right)$ we have a unique tuple $(\rho_i \mbox{ : } i \in Z_0)$, with $\rho_i \in {\rm Hom}_{\mathfrak{R}}([i], b)$. Then we have 
\begin{align*}
{\rm Card}(Z_1) & = \sum \left\{ {\rm Card}(Z^y_1) \mbox{ : } y \in {\rm Ker}\left( \partial_0 |_{\left< Z_0 \right>} \right) \right\} \\
& \leq {\rm Card}\left( {\rm Ker}\left( \partial_0 |_{\left< Z_0 \right>} \right) \right) \mbox{ since each $Z^y_1$ is finite}, \\
& \leq {\rm Card}\left( \{ (r_i \mbox{ : } i \in Z_0 ) \mbox{ : } y \in {\rm Ker}\left( \partial_0 |_{\left< Z_0 \right>} \right) \} \right) \\
& \leq \prod_{i \in Z_0} {\rm Card}\left( {\rm Hom}_{\mathfrak{R}} ([i], b) \right) \leq \kappa,
\end{align*}
and so $Z_1$ is a $\kappa$-small subset of $I_1$ such that $\partial_1 \left( \left< Z_1 \right> \right) \supseteq {\rm Ker}\left( \partial_0 |_{\left< Z_0 \right>} \right)$, i.e. $(\partial_1)_b \left( \left< Z_1 \right>(b) \right) \supseteq {\rm Ker}\left( (\partial_0)_b |_{\left< Z_0 \right>(b)} \right)$ for every $b \in {\rm Ob}(\mathfrak{R})$. In a similar way, we can find a $\kappa$-small subset $Z_2 \subseteq I_2$ such that $\partial_2\left( \left< Z_2 \right> \right) \supseteq {\rm Ker}\left( \partial_1 |_{\left< Z_1 \right>} \right)$. We keep repeating this procedure until we get a $\kappa$-small subset $Z_n \subseteq I_n$ such that $\partial \left( \left< Z_n \right> \right) \supseteq {\rm Ker}\left( \partial_{n-1} |_{\left< Z_{n-1} \right>} \right)$. 

The next step in the zig-zag procedure consists in choosing a $\kappa$-small subset $Z^{(1)}_{n-1} \subseteq I_{n-1}$, containing $Z_{n-1}$, such that $\partial_n \left( \left< Z_n \right> \right) \subseteq {< Z^{(1)}_{n-1}>}$. Let $y \in \partial_n \left( \left< Z_n \right> \right)$ of degree $b$. Then $y = (\partial_n)_b(z)$, where $z = \sum_{i \in Z_n} r_i \cdot [i]$. We have $y = \sum_{i \in Z_n} (\partial_n)_b(r_i \cdot [i])$. On the other hand, $(\partial_n)_b(r_i \cdot [i]) = \sum_{j \in Z^{(1), y, i}_{n-1}} q_j \cdot [j]$ for a finite $Z^{(1), y, i}_{n-1} \subseteq I_{n-1}$. Thus $y = \sum_{i \in Z_n} \sum_{j \in Z^{(1), y, i}_{n-1}} q_j \cdot [j] = \sum_{j \in Z^{(1)}_{n-1}} q_j \cdot [j]$, where $Z^{(1)}_{n-1}$ is the disjoint union $\bigsqcup \{ Z^{(1), y, i}_{n-1} \mbox{ : } y \in \partial_n \left( \left< Z_n \right> \right) \mbox{ and } i \in Z_n \}$. We have 
\begin{align*}
{\rm Card}\left(Z^{(1)}_{n-1}\right) & = \sum \left\{ {\rm Card}\left(Z^{(1), y, i}_{n-1}\right) \mbox{ : } y \in \partial_n \left( \left< Z_n \right> \right) \mbox{ and } i \in Z_n \right\} \\
& \leq {\rm Card}\left( \{ (r_i \mbox{ : } i \in Z_n) \mbox{ : } y \in \partial_n \left( \left< Z_n \right> \right) \} \right) \leq \kappa. 
\end{align*}
Then $Z^{(1)}_{n-1}$ is a $\kappa$-small subset of $I_{n-1}$ such that $\partial_n \left( \left< Z_n \right> \right) \subseteq {< Z^{(1)}_{n-1} >}$. Note that we can construct $Z^{(1)}_{n-1}$ containing $Z_{n-1}$. Similarly, there is a $\kappa$-small subset $Z^{(1)}_{n-2} \subseteq I_{n-2}$ containing $Z_{n-2}$ such that $\partial_{n-1} {(< Z^{(1)}_{n-1}>)} \subseteq {< Z^{(1)}_{n-2}>}$. 

At this point, we just need to mimic the argument given in the proof of \cite[Proposition 4.1]{Aldrich}, with the corresponding considerations for ${\rm {\bf Mod}}(\mathfrak{R})$. Set $J_k := Z_k \cup Z^{(1)}_{k} \cup \cdots$ for every $0 \leq k \leq n$. It is clear that $\left< J_k \right> := \bigoplus_{i \in J_k} {\rm Hom}_{\mathfrak{R}}([i], -)$ is a $\kappa$-small sub-module of $\left< I_k \right>$. By construction, we have an exact sequence \[ 0 \longrightarrow \left< J_n \right> \stackrel{\partial_n}\longrightarrow \left< J_{n-1} \right> \longrightarrow \cdots \longrightarrow \left< J_1 \right> \stackrel{\partial_1}\longrightarrow \left< J_0 \right> \stackrel{\partial_0}\longrightarrow N \longrightarrow 0, \] where $N = {\rm CoKer}(\left< J_1 \right>\stackrel{\partial_1}\longrightarrow \left< J_0 \right>)$. Note that $x \in N$ and that each $\left< J_k \right>$ is projective by Corollary \ref{projmodfrak}. It is only left to show that $M / N$ is $n$-projective. It suffices to take the quotient of the resolution of $M$ by the resolution of $N$, to get an exact sequence \[ 0 \longrightarrow \left< I_n\right> / \left<J_n\right> \longrightarrow \left<I_{n-1}\right> / \left<J_{n-1}\right> \longrightarrow \cdots \longrightarrow \left<I_0\right> / \left<J_0\right> \longrightarrow M / N \longrightarrow 0. \] It is not hard to check that $\left< I_k \right> / \left< J_k \right> \cong \left< I_k - J_k \right>$. So the previous sequence is a projective resolution of length $n$ of $M / N$.  
\end{proof}

\

\begin{proposition}\label{extensionespeques} Let $(\mathcal{A,B})$ be a hereditary cotorsion pair in ${\rm {\bf Mod}}(\mathfrak{R})$ and $\kappa$ be a regular cardinal satisfying $\kappa > {\rm Card}({\rm Hom}_{\mathfrak{R}}(a,b))$ for every $a, b \in {\rm Ob}(\mathfrak{R})$. Suppose for each left $n$-$\mathcal{A}$-module $X$ and each $x \in X$, there exists a $\kappa$-small left $n$-$\mathcal{A}$-module $X_x \subseteq X$ such that $x \in X_x$ and that $X / X_x$ is also a left $n$-$\mathcal{A}$-module. Then every left $n$-$\mathcal{A}$-module is a transfinite extension of the set of $\kappa$-small left $n$-$\mathcal{A}$-modules. 
\end{proposition}
\begin{proof} Let $X$ be a left $\mathfrak{R}$-module as described in the statement. Choose any $x_0 \in X$. Then there exists a small left $n$-$\mathcal{A}$-module $X_0$ such that $x_0 \in X_0$ and such that $X / X_0$ is also a left $n$-$\mathcal{A}$-module. Now choose a class $x_1 + X_0 \neq 0 + X_0$. Then there exists a small left $n$-$\mathcal{A}$-module $X_1 / X_0$ such that $x_1 + X_0 \in X_1/X_0$ and such that $X / X_1 \cong (X/X_0) / (X_1 / X_0)$ is left $n$-$\mathcal{A}$. Note that $X_0 \subseteq X_1$ and that $X_1$ is small since ${\rm Card}(X_1) = {\rm Card}(X_1 / X_0) \cdot {\rm Card}(X_0)$. Since we have a short exact sequence $0 \longrightarrow X_0 \longrightarrow X_1 \longrightarrow X_1 / X_0 \longrightarrow 0$, where $X_0$ and $X_1 / X_0$ are left $n$-$\mathcal{A}$-modules, by Lemma \ref{cerradoporexts}, $X_1$ is also a left $n$-$\mathcal{A}$-module. Using transfinite induction, we can construct a family of modules $(X_{\alpha} \mbox{ : } \alpha < \lambda)$, for some ordinal $\lambda$, such that $X$ is the transfinite extension $X = {\rm CoLim}_{\alpha < \lambda} X_\alpha$.   
\end{proof}

The following theorem is proven in \cite{Aldrich} for the category $\Modl$. The arguments carry over easily to the category of modules over ringoids. \\

\begin{theorem}[{\cite[Theorem 4.2]{Aldrich}} for modules over ringoids] \label{elparnprojectivo} The class of $n$-projective left $\mathfrak{R}$-modules $\mathcal{P}_n({\rm {\bf Mod}}(\mathfrak{R}))$ is the left half of a hereditary and complete cotorsion pair $\left(\mathcal{P}_n({\rm {\bf Mod}}(\mathfrak{R})), (\mathcal{P}_n({\rm {\bf Mod}}(\mathfrak{R})))^\perp\right)$ cogenerated by the set $(\mathcal{P}_n({\rm {\bf Mod}}(\mathfrak{R})))^{\leq \kappa}$.
\end{theorem}


\section{$n$-projective model structures}

The goal in this section is to construct a new Abelian model structure on $\Cadl$ where the weak equivalences are given by the quasi-isomorphisms and the trivial cofibrations by the monomorphisms with $n$-projective cokernels. The motivation of this problem comes from the case $n = 0$ studied by M. Hovey in \cite[Section 2.3]{Hovey}. 

We recall the following notions given by J. Gillespie in \cite[Definition 3.3]{Gillespie2}. Let $(\mathcal{A,B})$ be a cotorsion pair in an Abelian category $\mathcal{C}$. A chain complex $X$ in $\Complexes$ is:
\begin{itemize}
\item[{\bf (1)}] An \underline{$\mathcal{A}$-complex} if $X$ is exact and each cycle $Z_m(X)$ is in $\mathcal{A}$.

\item[{\bf (2)}] A \underline{dg-$\mathcal{A}$-complex} (\textquotedblleft dg\textquotedblright \ for \textquotedblleft differential graded\textquotedblright) if $X_m \in \mathcal{A}$ for each $m$ and every chain map $X \longrightarrow B$ is homotopic to zero whenever $B$ is a $\mathcal{B}$-complex. 

\item[{\bf (3)}] A \underline{dg-$\mathcal{B}$-complex} if $X_m \in \mathcal{B}$ for each $m$ and every chain map $A \longrightarrow X$ is homotopic to zero whenever $A$ is an $\mathcal{A}$-complex. 
\end{itemize}

Let $\widetilde{\mathcal{A}}$, ${\rm dg}\widetilde{\mathcal{A}}$ and ${\rm dg}\widetilde{\mathcal{B}}$ denote the class of $\mathcal{A}$-complexes, dg-$\mathcal{A}$-complexes and dg-$\mathcal{B}$-complexes, respectively. \\

\begin{example} By Proposition \ref{npicomplexes}, $\widetilde{\mathcal{P}_n(\mathcal{C})}$ and $\widetilde{\mathcal{I}_n(\mathcal{C})}$ are the classes of $n$-projective and $n$-injective complexes in $\Complexes$, respectively. \\
\end{example} 






Consider the class $\mathcal{P}_n$ of $n$-projective left $R$-modules. We know by \cite[Theorem 4.2]{Aldrich} that $(\mathcal{P}_n,(\mathcal{P}_n)\mbox{}^\perp)$ is a hereditary and complete cotorsion pair in $\Modl$. So by \cite[Corollary 3.8 and Theorem 3.12]{Gillespie2} we have two compatible cotorsion pairs $(\widetilde{\mathcal{P}_n},(\widetilde{\mathcal{P}_n})^\perp)$ and $({\rm dg}\widetilde{\mathcal{P}_n}, ({\rm dg}\widetilde{\mathcal{P}_n})^\perp)$. By Example \ref{particulares} {\bf (2)} and Proposition \ref{elparnprojectivo}, we have the following result. \\

\begin{corollary}\label{nprojpaircomp} The cotorsion pair $(\widetilde{\mathcal{P}_n},(\widetilde{\mathcal{P}_n})\mbox{}^{\perp})$ is hereditary and complete. 
\end{corollary}

The fact that the pair $({\rm dg}\widetilde{\mathcal{P}_n}, ({\rm dg}\widetilde{\mathcal{P}_n})^\perp)$ is complete shall be a consequence of the following result. \\

\begin{proposition}\label{compacogen} Let $(\mathcal{A} \cap \mathcal{W}, \mathcal{B})$ and $(\mathcal{A}, \mathcal{B} \cap \mathcal{W})$ be two compatible cotorsion pairs in an Abelian category $\mathcal{C}$.
\begin{itemize}
\item[{\bf (1)}] Suppose $(\mbox{}^\perp\mathcal{W},\mathcal{W})$ is a cotorsion pair cogenerated by a set $\mathcal{S}_{\mathcal{W}}$. If the pair $(\mathcal{A} \cap \mathcal{W}, \mathcal{B})$ is also cogenerated by a set $\mathcal{S}_{\mathcal{A} \cap \mathcal{W}}$, then $(\mathcal{A}, \mathcal{B} \cap \mathcal{W})$ is cogenerated by $\mathcal{S} = \mathcal{S}_{\mathcal{A} \cap \mathcal{W}} \cup \mathcal{S}_{\mathcal{W}}$.

\item[{\bf (2)}] Suppose $\mathcal{C}$ has enough projective and injective objects. If $(\mbox{}^\perp\mathcal{W},\mathcal{W})$ and $(\mathcal{W},\mathcal{W}^\perp)$ are complete cotorsion pairs, then $(\mathcal{A}, \mathcal{B} \cap \mathcal{W})$ is complete if, and only if, $(\mathcal{A} \cap \mathcal{W}, \mathcal{B})$ is.
\end{itemize}
\end{proposition}
\begin{proof} Part {\bf (1)} follows by $\mathcal{B} \cap \mathcal{W} = (\mathcal{S}_{\mathcal{A} \cap \mathcal{W}})^\perp \cap (\mathcal{S}_{\mathcal{W}})^\perp = (\mathcal{S}_{\mathcal{A} \cap \mathcal{W}} \cup \mathcal{S}_{\mathcal{W}})^\perp$. 

For Part {\bf (2)}, we only prove the implication $(\Longrightarrow)$, since the other is dual. So suppose $(\mathcal{A}, \mathcal{B} \cap \mathcal{W})$ is complete and let $X$ be an object in $\mathcal{C}$. Since $(\mathcal{W}, \mathcal{W}^\perp)$ is complete, there exists a short exact sequence $0 \longrightarrow X \longrightarrow C \longrightarrow W \longrightarrow 0$, where $C \in \mathcal{W}^\perp$ and $W \in \mathcal{W}$. Since $(\mathcal{A}, \mathcal{B} \cap \mathcal{W})$ is complete, there exists a short exact sequence $0 \longrightarrow B \longrightarrow A \longrightarrow W \longrightarrow 0$, where $A \in \mathcal{A}$ and $B \in \mathcal{B} \cap \mathcal{W}$. Taking the pullback of $C \longrightarrow W$ and $A \longrightarrow W$, we get two exact sequences $0 \longrightarrow X \longrightarrow C \times_W A \longrightarrow A \longrightarrow 0$ and $0 \longrightarrow B \longrightarrow C \times_W A \longrightarrow C \longrightarrow 0$ (See \cite[Proposition 2, page 203]{MacLane}). Since $\mathcal{W}$ is closed under extensions, we have $A \in \mathcal{A} \cap \mathcal{W}$. It suffices to show $C \times_W A \in \mathcal{B}$. Note that $\mathcal{A} \cap \mathcal{W} \subseteq \mathcal{W}$ implies $\mathcal{W}^\perp \subseteq (\mathcal{A} \cap \mathcal{W})^\perp = \mathcal{B}$. Then $C \in \mathcal{B}$. It follows $C \times_W A \in \mathcal{B}$ since $\mathcal{B}$ is closed under extensions. Hence $0 \longrightarrow X \longrightarrow C \times_W A \longrightarrow A \longrightarrow 0$ is a short exact sequence with $C \times_W A \in \mathcal{B}$ and $A \in \mathcal{A} \cap \mathcal{W}$. By \cite[Corollary 2.4]{Salce}, $(\mathcal{A} \cap \mathcal{W}, \mathcal{B})$ is complete.  
\end{proof}

\

\begin{remark} A cotorsion pair in an Abelian category $\mathcal{C}$ is \underline{functorially complete} if for every object $X$ there exists a special $\mathcal{A}$-pre-cover $A \longrightarrow X$ and a special $\mathcal{B}$-pre-envelope $X \longrightarrow B$ such that $A$ and $B$ are both functorial in $X$. 
\begin{itemize}
\item[{\bf (a)}] The assertion {\bf (2)} of the previous proposition is also valid for functorially complete cotorsion pairs. In this sense, {\bf (2)} can also be stated as follows:
\begin{itemize}
\item[{\bf (2')}] If $\mathcal{C}$ has functorially enough projective and injective objects, and if $(\mbox{}^\perp\mathcal{W},\mathcal{W})$ and $(\mathcal{W},\mathcal{W}^\perp)$ are functorially complete, then $(\mathcal{A}, \mathcal{B} \cap \mathcal{W})$ is functorially complete if, and only if, $(\mathcal{A} \cap \mathcal{W}, \mathcal{B})$ is. 
\end{itemize}

\item[{\bf (b)}] In the statement of the Eklof and Trlifaj Theorem, cotorsion pairs cogenerated by sets are actually functorially complete.

\item[{\bf (c)}] To prove the "if" part of the Hovey's correspondence, the pairs $(\mathcal{A} \cap \mathcal{W}, \mathcal{B})$ and $(\mathcal{A}, \mathcal{B} \cap \mathcal{W})$ need to be functorially complete (See the proof in \cite[Theorem 2.2]{Hovey2}.). \\
\end{itemize}
\end{remark}

To deduce that $({\rm dg}\widetilde{\mathcal{P}_n}, ({\rm dg}\widetilde{\mathcal{P}_n})^\perp)$ is complete, consider the class $\mathcal{E}$ of exact chain complexes. The trivial cotorsion pairs $(\mathcal{P}_0, \Modl)$ and $(\Modl, \mathcal{I}_0)$ induce two cotorsion pairs $({\rm dg}\widetilde{\mathcal{P}_0}, ({\rm dg}\widetilde{\mathcal{P}_0})^\perp)$ and $(\mbox{}^\perp({\rm dg}\widetilde{\mathcal{I}_0}), {\rm dg}\widetilde{\mathcal{I}_0})$ by \cite[Corollary 3.8]{Gillespie2}. In \cite[Propositions 2.3.4 \& 2.3.5]{Garcia}, it is proven that $({\rm dg}\widetilde{\mathcal{P}_0})^\perp = \mathcal{E} = \mbox{}^\perp({\rm dg}\widetilde{\mathcal{I}_0})$. On the one hand, $(\mathcal{E}, {\rm dg}\widetilde{\mathcal{I}_0})$ is complete by \cite[Theorem 2.3.17]{Garcia}. On the other hand, in \cite[Lemma 5.1]{Rada} it is proven that $({\rm dg}\widetilde{\mathcal{P}_0}, \mathcal{E})$ is cogenerated by a set, and so complete. Then putting $\mathcal{W} = \mathcal{E}$ in the previous proposition, we have that $({\rm dg}\widetilde{\mathcal{P}_n}, ({\rm dg}\widetilde{\mathcal{P}_n})^\perp)$ is a complete cotorsion pair (cogenerated by a set). Therefore, the following theorem follows by the Hovey's correspondence. \\

\begin{theorem}[The \underline{$n$-projective model structure on $\Cadl$}] \label{nprojmodeloncadl} There is a unique Abelian model structure on $\Cadl$ where the (trivial) cofibrations are the monomorphisms with cokernels in ${\rm dg}\widetilde{\mathcal{P}_n}$ (in $\widetilde{\mathcal{P}_n}$), the (trivial) fibrations are the epimorphisms with kernels in $( \widetilde{\mathcal{P}_n} )\mbox{}^\perp$ (in $({\rm dg}\widetilde{\mathcal{P}_n})\mbox{}^\perp$), and the trivial objects are the exact chain complexes. \\
\end{theorem}

\begin{remark} It is not hard to see that in the previous model structure, weak equivalences are given by quasi-isomorphisms. Moreover, in every Abelian model category with $\mathcal{E}$ as the class of trivial objects, a map is a weak equivalence if, and only if, it is a quasi-isomorphism. It is important to know that in \cite{Hovey2}, weak equivalences in Abelian model structures are defined as the composition of trivial cofibrations followed by trivial fibrations. 
\end{remark}

\newpage

An interesting question about a new model structure is whether it is monoidal or not. In order to have a monoidal model structure on $\mathcal{C}$, we need $\mathcal{C}$ to be equipped with a symmetric monoidal structure. Roughly speaking, a \underline{symmetric monoidal} \underline{structure} on $\mathcal{C}$ is given by a tensor product $- \otimes - : \mathcal{C} \times \mathcal{C} \longrightarrow \mathcal{C}$ and a unit object $S \in {\rm Ob}(\mathcal{C})$ such that for every $X, Y, Z \in {\rm Ob}(\mathcal{C})$, we have isomorphisms $X \otimes (Y \otimes Z) \cong (X \otimes Y) \otimes Z$, $X \otimes Y \cong Y \otimes X$, $S \otimes X \cong X$ and $X \otimes S \cong X$ (See the book \cite[Section 4.1]{Hovey} for details). The category $\Modl$, with $R$ a commutative ring, is an example of a symmetric monoidal category. The tensor product is the standard tensor product $\otimes_R$ of modules, where $R$ is the unit. 

We present two more examples of symmetric monoidal structures on $\Cadl$. Given two chain complexes $X \in {\rm Ob}(\Cadr)$ and $Y \in {\rm Ob}(\Cadl)$, the \underline{standard tensor product} of $X$ and $Y$ is given by the chain complex $X \otimes Y$ in ${\rm {\bf Ch}}({\rm {\bf Ab}})$ defined by $(X \otimes Y)_n := \bigoplus_{k \in \mathbb{Z}} X_k \otimes_R Y_{n-k}$, where the boundary maps are given by $x \otimes y \mapsto \partial^X_k(x) \otimes y + (-1)^k x \otimes \partial^Y_{n-k}(y)$, for every $x \otimes y \in X_k \otimes Y_{n-k}$. This construction defines a functor $- \otimes -$, from which one constructs the left derived functors ${\rm Tor}_i(-,-)$, with $i \geq 0$. From $\otimes$, the \underline{bar tensor product} of $X$ and $Y$ is defined as the chain complex $X \overline{\otimes} Y$ in ${\rm {\bf Ch}}({\rm {\bf Ab}})$ given by $(X \overline{\otimes} Y)_n := \frac{(X \otimes Y)_n}{B_n(X \otimes Y)}$, where the boundary maps are defined by $x \otimes y + B_n(X \otimes Y) \mapsto \partial^X_k(x) \otimes y + B_{n-1}(X \otimes Y)$, for every $x \otimes y + B_n(X \otimes Y) \in (X \overline{\otimes} Y)_n$. As far as the author knows, the definition of this tensor product first appeared in \cite{Garcia}. If $R$ is a commutative ring, then $(\Cadl, \otimes)$ and $(\Cadl, \overline{\otimes})$ are symmetric monoidal categories, where $S^0(R)$ and $D^1(R)$ are the units with respect to $\otimes$ and $\overline{\otimes}$, respectively (See \cite[Proposition 4.2.13]{Hovey} and \cite[Proposition 4.2.1 4]{Garcia}.). 

A \underline{monoidal model category} is a model category $\mathcal{C}$ equipped with a symmetric monoidal structure $(\otimes, S)$ such that the following conditions are satisfied: 
\begin{itemize}
\item[{\bf (1)}] For every pair of maps $f : U \longrightarrow V$ and $g : W \longrightarrow X$, the pushout of $f \otimes {\rm id}_W$ and ${\rm id}_U \otimes g$ induces a map $f \square g : (V \otimes W) \coprod_{U \otimes W} (U \otimes X) \rightarrow V \otimes X$ making the following diagram commute: \\
\[ \begin{tikzpicture}
\matrix (m) [matrix of math nodes, row sep=2em, column sep=2em]
{ U \otimes W & U \otimes X \\ V \otimes W & (V \otimes W) \coprod_{U \otimes W} (U \otimes X) \\ & & V \otimes X \\  };
\path[->]
(m-1-1) edge node[above] {${\rm id}_U \otimes g$} (m-1-2) edge node[left] {$f \otimes {\rm id}_W$} (m-2-1) (m-1-2) edge (m-2-2) (m-2-1) edge (m-2-2)
(m-1-2) edge [bend left=30] (m-3-3) (m-2-1) edge [bend right=30] (m-3-3);
\path[dotted,->]
(m-2-2) edge node[sloped] {$f \square g$} (m-3-3);
\end{tikzpicture} \]
If, given cofibrations $f : U \rightarrow V$ and $g : W \rightarrow X$ in $\mathcal{C}$, the induced map $f \square g$ is a cofibration, which is trivial if either $f$ or $g$ is. \\

\item[{\bf (2)}] Using functorial factorizations, write $0 \rightarrow S = 0 \longrightarrow Q(S) \stackrel{q}\longrightarrow S$ as the composition of a cofibration followed by a trivial fibration. Then the maps $q \otimes X : Q(S) \otimes X \longrightarrow S \otimes X$ and $X \otimes q : X \otimes Q(S) \longrightarrow X \otimes S$ are weak equivalences for all cofibrant objects $X$. 
\end{itemize}

It is known that the $n$-projective model structure is monoidal with respect to the usual tensor product on $\Cadl$ if $n = 0$ (See \cite[Chapter 4]{Hovey}). In general, this is not the case if $n > 0$. For instance, consider $n = 1$ and $R = \mathbb{Z}$. Note that $\mathbb{Z}_2$ is a $1$-projective $\mathbb{Z}$-module, since it is not projective and there exists a short exact sequence $0 \longrightarrow \mathbb{Z} \stackrel{2 \times}\longrightarrow \mathbb{Z} \stackrel{\pi}\longrightarrow \mathbb{Z}_2 \longrightarrow 0$ where $2 \times$ is the map $x \mapsto 2 \times x$ and $\pi$ is the canonical projection $x \mapsto \overline{x} \in \{ 0, 1 \}$.  Now let $X$ be the complex given by the previous sequence, where $X_1 = \mathbb{Z}$, $X_0 = \mathbb{Z}$ and $X_{-1} = \mathbb{Z}_2$. We have $X \in \widetilde{\mathcal{P}_1}$. Consider also the complex $S^0(\mathbb{Z}_2)$. Then $S^0(\mathbb{Z}_2) \in {\rm dg}\widetilde{\mathcal{P}_1}$ since ${\rm Ext}^1_{{\rm {\bf Ch}}({\rm {\bf Mod}}_{\mathbb{Z}})}(S^0(\mathbb{Z}_2), Y) \cong {\rm Ext}^1_{\mathbb{Z}}(\mathbb{Z}_2, Z_0 Y)$ (by \cite[Lemma 4.2]{Gillespie}) and ${\rm Ext}^1_{\mathbb{Z}}(\mathbb{Z}_2, Z_0 Y) = 0$ for every $Y \in \widetilde{\mathcal{P}_1\mbox{}^\perp} = ({\rm dg}\widetilde{\mathcal{P}_1})\mbox{}^\perp$. It is not hard to see that
\begin{align*}
(S^0(\mathbb{Z}_2) \otimes X)_m & = \left\{ \begin{array}{ll} \mathbb{Z}_2 \otimes_{\mathbb{Z}} \mathbb{Z} & \mbox{ if }m = 1, \\ \mathbb{Z}_2 \otimes_{\mathbb{Z}} \mathbb{Z} & \mbox{ if }m =0, \\ \mathbb{Z}_2 \otimes_{\mathbb{Z}} \mathbb{Z}_2 & \mbox{ if }m = -1, \\ 0 & \mbox{ otherwise}. \end{array} \right.
\end{align*}
and that $\partial^{S^0(\mathbb{Z}_2) \otimes X}_1$ is the zero map, so the sequence \[ \cdots \longrightarrow 0 \longrightarrow \mathbb{Z}_2 \otimes_{\mathbb{Z}} \mathbb{Z} \stackrel{\partial_1^{S^0(\mathbb{Z}_2) \otimes X}}\longrightarrow \mathbb{Z}_2 \otimes_{\mathbb{Z}} \mathbb{Z} \stackrel{\partial_0^{S^0(\mathbb{Z}_2) \otimes X}}\longrightarrow \mathbb{Z}_2 \otimes_{\mathbb{Z}} \mathbb{Z}_2 \longrightarrow 0 \longrightarrow \cdots \] is not exact. We have a cofibration $0 \longrightarrow S^0(\mathbb{Z}_2)$ and a trivial cofibration $0 \longrightarrow X$, but the induced map $(0 \longrightarrow S^0(\mathbb{Z}_2)) \square (0 \longrightarrow X) = 0 \longrightarrow S^0(\mathbb{Z}_2) \otimes X$ is cofibration but not a weak equivalence, since $S^0(\mathbb{Z}_2) \otimes X$ is not exact. Therefore, the $n$-projective model structure on ${\rm {\bf Ch}}({\rm {\bf Mod}}_{\mathbb{Z}})$ is not monoidal with respect to the tensor product $\otimes$. We can conclude the same for ${\rm {\bf Ch}}({\rm {\bf Mod}}_{\mathbb{Z}})$ with respect to the tensor product $\overline{\otimes}$, since 
\begin{align*}
(S^0(\mathbb{Z}_2) \overline{\otimes} X)_m & = \left\{ \begin{array}{ll} \mathbb{Z}_2 \otimes_{\mathbb{Z}} \mathbb{Z} & \mbox{ if $m = 1$}, \\ \mathbb{Z}_2 \otimes_{\mathbb{Z}} \mathbb{Z} & \mbox{ if $m = 0$}, \\ 0 & \mbox{ otherwise}. \end{array} \right. 
\end{align*}
and $S^0(\mathbb{Z}_2) \overline{\otimes} X = \cdots \longrightarrow 0 \longrightarrow \mathbb{Z}_2 \otimes_{\mathbb{Z}} \mathbb{Z} \stackrel{\partial_1^{S^0(\mathbb{Z}_2) \overline{\otimes} X}}\longrightarrow \mathbb{Z}_2 \otimes_{\mathbb{Z}} \mathbb{Z} \longrightarrow 0 \longrightarrow \cdots$, where $\partial^{S^0(\mathbb{Z}_2) \overline{\otimes} X}_1$ is the zero map. Hence $S^0(\mathbb{Z}_2) \overline{\otimes} X$ is not an exact complex, and so the induced chain map $0\longrightarrow S^0(\mathbb{Z}_2) \overline{\otimes} X = (0 \longrightarrow S^0(\mathbb{Z}_2)) \square (0 \longrightarrow X)$ is not a trivial cofibration. 

We conclude this section presenting the dual of the $n$-projective model structure. The case of the injective dimension is easier to study. Assume for the rest of this section that $\mathcal{C}$ is a Grothendieck category with a generator $G$. Grothendieck categories have enough injective objects. This was proven by A. Grothendieck in his famous Tohoku paper {\it Sur quelques points d'alg\`ebre homologique} (See \cite[Theorem 1.10.1]{Grothendieck}.). Moreover, in this paper it is proven the validity of the Baer Criterion for any Grothendieck category (See \cite[Lemme 1, page 136]{Grothendieck}.). The following result follows. \\

\begin{proposition} Let $\mathcal{C}$ be a Grothendieck category with a generator $G$. Then an object $Y$ of $\mathcal{C}$ is $n$-injective if, and only if, ${\rm Ext}^{n+1}_{\mathcal{C}}(G/J,Y) = 0$ for every subobject $J$ of $G$. 
\end{proposition}

The following follows easily. \\

\begin{corollary}\label{theninjpair} If $\mathcal{C}$ is a Grothendieck category with generator $G$, then $\mathcal{I}_n(\mathcal{C})$ is the right half of a cotorsion pair $(\mbox{}^\perp(\mathcal{I}_n(\mathcal{C})), \mathcal{I}_n(\mathcal{C}))$ cogenerated by the set of all $S \in \Omega^n(G/J)$ with $J$ running over the set of all subobjects of $G$. \\
\end{corollary}

It is not hard to see that the class $\mathcal{I}_n(\mathcal{C})$ is coresolving, and hence the pair $(\mbox{}^\perp(\mathcal{I}_n(\mathcal{C})), \mathcal{I}_n(\mathcal{C}))$ is also hereditary. If in addition $\mathcal{C}$ has enough enough projective objects, \cite[Corollary 3.8]{Gillespie2} and Proposition \ref{theninjpair} imply that $({\rm dg}\widetilde{\mbox{}^\perp(\mathcal{I}_n(\mathcal{C}))}, \widetilde{\mathcal{I}_n(\mathcal{C})})$ and $(\widetilde{\mbox{}^\perp(\mathcal{I}_n(\mathcal{C}))}, {\rm dg}\widetilde{\mathcal{I}_n(\mathcal{C})})$ are compatible cotorsion pairs in $\Complexes$. Moreover, by Proposition \ref{npicomplexes}, we know $\widetilde{\mathcal{I}_n(\mathcal{C})} = \mathcal{I}_n(\Complexes)$. It follows by Corollary \ref{theninjpair}, the pair $({\rm dg}\widetilde{\mbox{}^\perp(\mathcal{I}_n(\mathcal{C}))}, \widetilde{\mathcal{I}_n(\mathcal{C})})$ is complete. Hence, by Proposition \ref{compacogen}, $(\widetilde{\mbox{}^\perp(\mathcal{I}_n)}, {\rm dg}\widetilde{\mathcal{I}_n})$ is a complete cotorsion pair. Then the following theorem follows. 

\vspace{0.5cm}

\begin{theorem}[The \underline{$n$-injective model structure on $\Cadl$}] \label{ninjmodel} There exists a unique Abelian model structure on $\Cadl$, where the (trivial) fibrations are the epimorphisms with kernel in ${\rm dg}\widetilde{\mathcal{I}_n}$ (in $\widetilde{\mathcal{I}_n}$), the (trivial) cofibrations are the monomorphisms with cokernel in $\mbox{}^\perp(\widetilde{\mathcal{I}_n})$ (in $\mbox{}^\perp({\rm dg}\widetilde{\mathcal{I}_n})$), and the weak equivalences are the quasi-isomorphisms. 
\end{theorem}


\section{Degreewise $n$-projective model structures}

This section is devoted to construct an Abelian model structure on $\Cadl$ where the exact complexes are the trivial objects, and the cofibrant objects are the chain complexes whose terms have projective dimension at most $n$. This structure represents an extension of the model structure found by E. E. Enochs and coauthors in \cite[Theorem 5.5]{Rada}. For this purpose, it is useful to recall more induced cotorsion pairs in $\Complexes$ from a cotorsion pair $(\mathcal{A,B})$ in $\mathcal{C}$, along with some of their properties.

Let $\mathcal{A}$ be a class of objects in an Abelian category $\mathcal{C}$. Let ${\rm dw}\widetilde{\mathcal{A}}$ denote the class of all complexes $X$ in $\Complexes$ such that $X_n \in \mathcal{A}$. We shall say that a complex in ${\rm dw}\widetilde{\mathcal{A}}$ is a \underline{degreewise $\mathcal{A}$-complex}, or a \underline{dw-$\mathcal{A}$-complex}. The class ${\rm ex}\widetilde{\mathcal{A}} := {\rm dw}\widetilde{\mathcal{A}} \cap \mathcal{E}$ shall be called the class of \underline{exact dw-$\mathcal{A}$-complexes}. This definition was given by J. Gillespie in \cite[Definition 3.1]{Gillespie}. 

In \cite[Proposition 3.2]{Gillespie}, it is proven by Gillespie that if $(\mathcal{A,B})$ be a cotorsion pair in an Abelian category $\mathcal{C}$, then $({\rm dw}\widetilde{\mathcal{A}}, ({\rm dw}\widetilde{\mathcal{A}})^\perp)$ and $(\mbox{}^\perp({\rm dw}\widetilde{\mathcal{B}}), {\rm dw}\widetilde{\mathcal{B}})$ are cotorsion pairs in $\Complexes$. The same author proved in \cite[Proposition 3.3]{Gillespie} that if in addition $\mathcal{B}$ contains a cogenerator of finite injective dimension, then $({\rm ex}\widetilde{\mathcal{A}}, ({\rm ex}\widetilde{\mathcal{A}})^\perp)$ is a cotorsion pair. This is also true if we assume instead that $\mathcal{C}$ has enough injective objects. Dually, if either $\mathcal{A}$ contains a generator of finite projective dimension or $\mathcal{C}$ has enough projective objects, then $(\mbox{}^\perp({\rm ex}\widetilde{\mathcal{B}}), {\rm ex}\widetilde{\mathcal{B}})$ is a cotorsion pair in $\Complexes$. 

Considering the trivial cotorsion pair $(\mathcal{P}_0, \Modl)$ in $\Modl$, we have two cotorsion pairs $({\rm dw}\widetilde{\mathcal{P}_0}, ({\rm dw}\widetilde{\mathcal{P}_0})\mbox{}^\perp)$ and $({\rm ex}\widetilde{\mathcal{P}_0}, ({\rm ex}\widetilde{\mathcal{P}_0})\mbox{}^\perp)$ in $\Cadl$. In \cite{Rada}, the authors show that these two pairs are cogenerated by sets. Concerning the former pair, they prove in \cite[Theorem 4.4]{Rada} that if $P \in {\rm dw}\widetilde{\mathcal{P}_0}$ and if $\mathcal{S}$ is a set of representatives of complexes with all terms countably generated projective modules, then $P$ has an $\mathcal{S}$-filtration. For the latter pair, if $\kappa$ is an infinite cardinal with $\kappa > {\rm Card}(R)$, then every $P \in {\rm ex}\widetilde{\mathcal{P}_0}$ has a $({\rm ex}\widetilde{\mathcal{P}_0})^{\leq \kappa}$-filtration (See \cite[Theorem 4.6]{Rada}.). It follows by Proposition \ref{gentransex} that $({\rm dw}\widetilde{\mathcal{P}_0}, ({\rm dw}\widetilde{\mathcal{P}_0})\mbox{}^\perp)$ and $({\rm ex}\widetilde{\mathcal{P}_0}, ({\rm ex}\widetilde{\mathcal{P}_0})\mbox{}^\perp)$ are complete cotorsion pairs. In \cite{Rada} it is also proven that these pairs are compatible. So it follows by the Hovey's correspondence the existence of a model structure on $\Cadl$, which shall be referred as the \underline{degreewise projective model structure}, where the (trivial) cofibrations are the monomorphisms with cokernels in ${\rm dw}\widetilde{\mathcal{P}_0}$ (in ${\rm ex}\widetilde{\mathcal{P}_0}$), the (trivial) fibrations are the epimorphisms with kernels in $({\rm ex}\widetilde{\mathcal{P}_0})\mbox{}^\perp$ (in $({\rm dw}\widetilde{\mathcal{P}_0})\mbox{}^\perp$), and the weak equivalences are the quasi-isomorphisms (See {\cite[Theorem 5.5]{Rada}}). 

We generalize the arguments given in \cite{Rada} to any projective dimension. In other words, we shall prove that every (exact) degreewise $n$-projective complex is filtered by some set, in order to obtain the following model structure: \\

\begin{theorem}[The \underline{degreewise $n$-projective model structure}] \label{nprojgrad} There exists a unique Abelian model structure on $\Cadl$ where the weak equivalences are the quasi-isomorphisms, the (trivial) cofibrations are the monomorphisms with cokernels in ${\rm dw}\widetilde{\mathcal{P}_n}$ (in ${\rm ex}\widetilde{\mathcal{P}_n}$), and the (trivial) fibrations are the epimorphisms with kernels in $({\rm ex}\widetilde{\mathcal{P}_n})^\perp$ (in $({\rm dw}\widetilde{\mathcal{P}_n})^\perp$). \\
\end{theorem}

 We first need to note that by \cite[Proposition 3.2]{Gillespie}, we have two cotorsion pairs $({\rm dw}\widetilde{\mathcal{P}_n}, ({\rm dw}\widetilde{\mathcal{P}_n})\mbox{}^\perp)$ and $({\rm ex}\widetilde{\mathcal{P}_n}, ({\rm ex}\widetilde{\mathcal{P}_n})\mbox{}^\perp)$. \\
 
\begin{definition} We say that a chain complex $X$ in $\Cadl$ is \underline{degreewise} \underline{$n$-projective} if $X \in {\rm dw}\widetilde{\mathcal{P}_n}$. \\
\end{definition}

The proof of \cite[Theorem 4.4]{Rada} is based in the following result by I. Kaplansky on projective modules. Recall that a left $R$-module is said to be \underline{countably generated} if it is generated as a module by a countable subset. \\

\begin{theorem}[I. Kaplansky. See \cite{Kaplansky}] If $P$ is a projective module then $P$ is the direct sum of countably generated projective modules. \\
\end{theorem} 

So when one thinks of a possible generalization of the degreewise projective model structure for $n$-projective modules, a good question would be if it is possible to generalize the Kaplansky's Theorem for such modules. 

Let $M \in \mathcal{P}_n$ be an $n$-projective module, so we can choose an projective resolution of length $n$, say $0 \longrightarrow P_n \longrightarrow \cdots \longrightarrow P_1 \longrightarrow P_0 \longrightarrow M \longrightarrow 0$. By Kaplansky's Theorem we can write $P_k = \bigoplus_{i \in I_k} P^i_k$, where $P^i_k$ is a countably generated projective module, for every $i \in I_k$ and every $0 \leq k \leq n$. Then we can rewrite the previous resolution as \[ 0 \longrightarrow \bigoplus_{i \in I_n} P^i_n \longrightarrow \cdots \longrightarrow \bigoplus_{i \in I_1} P^i_1 \longrightarrow \bigoplus_{i \in I_0} P^i_0 \longrightarrow M \longrightarrow 0. \] From now on we shall write any projective resolution of length $n$ by using such direct sum decompositions. We shall denote by $\mathcal{P}_n^{\aleph_0}$ the set of all modules $M$ having a projective resolution as above, where $I_k$ is a countable set for each $0 \leq k \leq n$. 

We shall prove that $({\rm dw}\widetilde{\mathcal{P}_n}, ({\rm dw}\widetilde{\mathcal{P}_n})^\perp)$ is a cotorsion pair cogenerated by the set ${\rm dw}\widetilde{\mathcal{P}_n^{\aleph_0}}$, by showing that every complex in ${\rm dw}\widetilde{\mathcal{P}_n}$ is ${\rm dw}\widetilde{\mathcal{P}_n^{\aleph_0}}$-filtered. We need the following extension of  Kaplansky's Theorem for $n$-projective modules. \\

\begin{proposition}[Kaplansky's Theorem fon $n$-projective modules]\label{Kaplansky} Let $R$ be a Noetherian ring. Let $M \in \mathcal{P}_n$ and $N \subseteq M$ be a countably generated sub-module of $M$. Then there exists a $\mathcal{P}^{\aleph_0}_n$-filtration of $M$, $(M_\alpha : \alpha < \lambda)$  with $\lambda > 1$, such that $M_1 \in \mathcal{P}^{\aleph_0}_n$ and $N \subseteq M_1$. 
\end{proposition}
\begin{proof} Let $M \in \mathcal{P}_n$ with a finite projective resolution where each projective term $P_k$ is written as a direct sum $\bigoplus_{i \in I_k} P^i_k$ of countably generated projective modules, for every $0 \leq i \leq n$. We shall construct a $\mathcal{P}^{\aleph_0}_n$-filtration $(M_\alpha \mbox{ : } \alpha < \lambda)$ of $M$, with $N \subseteq M_1$, by using  the zig-zag procedure described in the proof of Lemma \ref{lemaoyogen}. 

Set $M_0 = 0$. To construct $M_1$, consider a countable set $\mathcal{G}$ of generators of $N$. Since $f_0$ is surjective, for every $g \in \mathcal{G}$ we can choose $y_g \in \bigoplus_{i \in I_0} P^i_0$ such that $g = f_0(y_g)$. Consider the set $Y = \{ y_g \mbox{ : } g \in \mathcal{G} \}$. Since $Y$ is a countable subset of $\bigoplus_{i \in I_0} P^i_0$, we have that $\left< Y \right>$ is a countably generated sub-module of $P_0$. Choose a countable subset $I_0^{1, 0} \subseteq I_0$ such that $\left< Y \right> \subseteq \bigoplus_{i \in I^{1, 0}_0} P^i_0$. Then $f_0\left( \left< Y \right> \right) \subseteq N$. Consider ${\rm Ker}( f_0|_{\bigoplus_{i \in I_0^{1, 0}} P^i_0})$. Since $\bigoplus_{i \in I_0^{1, 0}} P^i_0$ is countably generated and ${\rm Ker}(f_0|_{\bigoplus_{i \in I_0^{1, 0}} P^i_0})$ is a sub-module of $\bigoplus_{i \in I_0^{1, 0}} P^i_0$, we have that ${\rm Ker}(f_0|_{\bigoplus_{i \in I_0^{1, 0}} P^i_0})$ is also countably generated, since $R$ is Noetherian. Let $\mathcal{B}$ be a countable set of generators of ${\rm Ker}(f_0|_{\bigoplus_{i \in I_0^{1, 0}} P^i_0})$. Let $b \in \mathcal{B}$, then $f(b) = 0$ and by exactness of the above sequence there exists $y_b \in \bigoplus_{i \in I_1} P^i_1$ such that $b = f_1(y_b)$. Let $Y' = \{ y_b \mbox{ : } b \in \mathcal{B} \}$. Note that $Y'$ is a countable subset of $(f_1)^{-1}({\rm Ker}(f_0|_{\bigoplus_{i \in I_0^{1, 0}} P^i_0}))$. Then $\left< Y' \right>$ is a countably generated sub-module of $\bigoplus_{i \in I_1} P^i_1$. Hence there exists a countable subset $I_1^{1, 0} \subseteq I_1$ such that $\bigoplus_{i \in I_1^{1, 0}} P^i_1 \supseteq \left< Y' \right>$. Thus $f_1(\bigoplus_{i \in I^{1, 0}_1} P^i_1) \supseteq f_1(\left< Y' \right>) \supseteq {\rm Ker}(f_0|_{\bigoplus_{i \in I_0^{1, 0}} P^i_0})$. Repeat the same argument until find a countable subset $I_n^{1, 0} \subseteq I_n$ such that $f_n(\bigoplus_{i \in I_n^{1, 0}} P^i_n) \supseteq {\rm Ker}( f_{n-1}|_{\bigoplus_{i \in I_{n-1}^{1, 0}} P^i_{n-1}})$. Now, $f_n(\bigoplus_{i \in I_n^{1,0}} P^i_n)$ is a countably generated sub-module of $\bigoplus_{i \in I_{n-1}} P^i_{n-1}$. Then choose a countable subset $I_{n-1}^{1,0} \subseteq I_{n-1}^{1, 1} \subseteq I_{n-1}$ such that $f_n(\bigoplus_{i \in I_n^{1, 0}} P^i_n) \subseteq \bigoplus_{i \in I_{n-1}^{1,1}} P^i_{n-1}$. Continue this zig-zag procedure infinitely many times and set $I^1_k = \bigcup_{m \geq 0} I_k^{1, m}$, for every $0 \leq k \leq n$. By construction, we get an exact sequence \[ 0 \longrightarrow \bigoplus_{i\in I_n^1} P^i_n \longrightarrow \bigoplus_{i \in I_{n-1}^1} P^i_{n-1} \longrightarrow \cdots \longrightarrow \bigoplus_{i \in I^1_1} P^i_1 \longrightarrow \bigoplus_{i \in I^1_0} P^i_0 \longrightarrow M_1 \longrightarrow 0 \] where $x \in M_1 := {\rm CoKer}(\bigoplus_{i \in I^1_1} \longrightarrow \bigoplus_{i \in I_0^1} P^i_0) \subseteq M$ and $N \subseteq M_1$. We take the quotient of the resolution of $M$ by the resolution of $M'$, and get a sequence \[ 0 \longrightarrow \bigoplus_{i \in I_n - I^1_n} P^i_n \longrightarrow \cdots \longrightarrow \bigoplus_{i \in I_1 - I^1_1} P^i_1 \longrightarrow \bigoplus_{i \in I_0 - I^1_0} P^i_0 \longrightarrow \frac{M}{M_1} \longrightarrow 0, \] which is exact since the class of exact complexes is thick. Hence we have a projective resolution of length $n$ for $M / M_1$. The rest of the proof follows by using an argument like in the proof of Proposition \ref{extensionespeques}.
\end{proof}

\

\begin{lemma} If $R$ is a Noetherian ring, then every chain complex in ${\rm dw}\widetilde{\mathcal{P}_n}$ has a ${\rm dw}\widetilde{\mathcal{P}_n^{\aleph_0}}$-filtration. 
\end{lemma}
\begin{proof} Let $X = (\cdots \longrightarrow X_{k+1} \stackrel{\partial_{k+1}}\longrightarrow X_k \stackrel{\partial_k}\longrightarrow X_{k-1} \longrightarrow \cdots)$ be a complex in ${\rm dw}\widetilde{\mathcal{P}_n}$. For each $k$ one has a projective resolution of $X_k$ of length $n$: \[ 0 \longrightarrow \bigoplus_{i \in I_n(k)} P^i_n(k) \longrightarrow \cdots \longrightarrow \bigoplus_{i \in I_1(k)} P^i_1(k) \longrightarrow \bigoplus_{i \in I_0(k)} P^i_0(k) \longrightarrow X_k \longrightarrow 0, \] where the modules appearing in the direct sums are countably generated and projective. It suffices to construct a nonzero sub-complex $X' \subseteq X$ such that each $X'_k$ has a projective resolution \[ 0 \longrightarrow \bigoplus_{i \in I'_n(k)} P^i_n(k) \longrightarrow \cdots \longrightarrow \bigoplus_{i \in I'_1(k)} P^i_1(k) \longrightarrow \bigoplus_{i \in I'_0(k)} P^i_0(k) \longrightarrow X'_k \longrightarrow 0, \] where $I'_i(k)$ is countable for every $0 \leq i \leq n$, and such that $X / X' \in {\rm dw}\widetilde{\mathcal{P}_n}$. 

Fix $m \in \mathbb{Z}$. Let $S$ be a countably generated sub-module of $X_m$. By the previous lemma, there exists a sub-module $\mathcal{P}^{\aleph_0}_n \ni X'_m \subseteq X_m$ such that $S \subseteq X'_m$. Note that $X'_m$ is also countably generated. Then $\partial_m(X'_m)$ is a countably generated sub-module of $X_{m-1}$, and so there exists $\mathcal{P}^{\aleph_0}_n \ni X'_{m-1} \subseteq X_{m-1}$ such that $\partial_m(X'_m) \subseteq X'_{m-1}$. Repeat the same procedure infinitely many times in order to obtain a sub-complex $X' = \cdots \longrightarrow X'_{k+1} \longrightarrow X'_k \longrightarrow X'_{k-1} \longrightarrow \cdots$  of $X$ such that $X'_k \in \mathcal{P}^{\aleph_0}_n$ for every $k \in \mathbb{Z}$ (We are setting $X'_k = 0$ for every $k > m$.). Hence $X' \in {\rm dw}\widetilde{\mathcal{P}^{\aleph_0}_n}$. Note from the proof of the previous lemma that the quotient $X / X'$ is in ${\rm dw}\widetilde{\mathcal{P}_n}$. We have for every $k \leq m$ one has the following projective resolutions of length $n$ for $X_k / X'_{k}$:  \[ 0 \longrightarrow \bigoplus_{i \in I_n(k) - I'_n(k)} P^i_n(k) \longrightarrow \cdots \longrightarrow \bigoplus_{i \in I_0(k) - I'_0(k)} P^i_0(k) \longrightarrow \frac{X_k}{X'_k} \longrightarrow 0. \] The rest of the proof follows by transfinite induction.
\end{proof}

\

It follows by the previous lemma and Proposition \ref{gentransex} that $({\rm dw}\widetilde{\mathcal{P}_n}, ({\rm dw}\widetilde{\mathcal{P}_n})^\perp)$ is a complete cotorsion pair in $\Cadl$ cogenerated by the set ${\rm dw}\widetilde{\mathcal{P}_n^{\aleph_0}}$, provided $R$ is a Noetherian ring. Under this hypothesis, it is easy to show that the pair $({\rm ex}\widetilde{\mathcal{P}_n}, ({\rm ex}\widetilde{\mathcal{P}_n})^\perp)$ is also complete using Proposition \ref{compacogen}, as long as we show that it is compatible with $({\rm dw}\widetilde{\mathcal{P}_n}, ({\rm dw}\widetilde{\mathcal{P}_n})^\perp)$. 

We know by \cite[Theorem 3.12]{Gillespie2} that the cotorsion pairs $( {\rm dg}\widetilde{\mathcal{A}}, \widetilde{\mathcal{B}} )$ and $( \widetilde{\mathcal{A}}, {\rm dg}\widetilde{\mathcal{B}} )$ are compatible if the inducing pair $(\mathcal{A},\mathcal{B})$ is hereditary. The author does not know if the same holds for the cotorsion pairs $( {\rm dw}\widetilde{\mathcal{A}}, ( {\rm dw}\widetilde{\mathcal{A}} )\mbox{}^\perp )$ and $( {\rm ex}\widetilde{\mathcal{A}}, ( {\rm ex}\widetilde{\mathcal{A}} )\mbox{}^\perp )$. 

Consider the case $\mathcal{C} = {\Modl}$. Since ${\rm dg}\widetilde{\mathcal{P}_0(\mathcal{C})} \subseteq {\rm dw}\widetilde{\mathcal{A}}$ and ${\rm ex}\widetilde{\mathcal{A}} \subseteq {\rm dw}\widetilde{\mathcal{A}}$, we have $({\rm dw}\widetilde{\mathcal{A}})\mbox{}^\perp \subseteq ( {\rm dg}\widetilde{\mathcal{P}_0(\mathcal{C})} )\mbox{}^\perp = \mathcal{E}$ and $( {\rm dw}\widetilde{\mathcal{A}} )\mbox{}^\perp \subseteq ( {\rm ex}\widetilde{\mathcal{A}} )\mbox{}^\perp$, and so $( {\rm dw}\widetilde{\mathcal{A}} )\mbox{}^\perp \subseteq ( {\rm ex}\widetilde{\mathcal{A}} )\mbox{}^\perp \cap \mathcal{E}$. Similarly, $\mbox{}^\perp( {\rm dw}\widetilde{\mathcal{B}} ) \subseteq \mbox{}^\perp( {\rm ex}\widetilde{\mathcal{B}} ) \cap \mathcal{E}$. The next result provides a relationship between completeness of these degreewise cotorsion pairs and the remaining inclusions. \\

\begin{proposition}\label{degreecomp} Let $(\mathcal{A,B})$ be a cotorsion pair in $\Cadl$. 
\begin{itemize}
\item[{\bf (1)}] If the pair $({\rm dw}\widetilde{\mathcal{A}}, ({\rm dw}\widetilde{\mathcal{A}})\mbox{}^\perp)$ is complete, then $({\rm dw}\widetilde{\mathcal{A}})\mbox{}^\perp = ({\rm ex}\widetilde{\mathcal{A}})\mbox{}^\perp \cap \mathcal{E}$.

\item[{\bf (2)}] If the pair $(\mbox{}^\perp({\rm dw}\widetilde{\mathcal{B}}), {\rm dw}\widetilde{\mathcal{B}})$ is complete, then $\mbox{}^\perp({\rm dw}\widetilde{\mathcal{B}}) = \mbox{}^\perp({\rm ex}\widetilde{\mathcal{B}}) \cap \mathcal{E}$.  
\end{itemize}
\end{proposition}
\begin{proof}
We only prove the first statement. It suffices to show the inclusion $( {\rm ex}\widetilde{\mathcal{A}} )\mbox{}^\perp \cap \mathcal{E} \subseteq ( {\rm dw}\widetilde{\mathcal{A}})\mbox{}^\perp$. Let $Y \in ( {\rm ex}\widetilde{\mathcal{A}} )\mbox{}^\perp \cap \mathcal{E}$. Since $( {\rm dw}\widetilde{\mathcal{A}}, ( {\rm dw}\widetilde{\mathcal{A}} )\mbox{}^\perp )$ is complete, there is a short exact sequence $0 \longrightarrow Y \longrightarrow X \longrightarrow A \longrightarrow 0$, with $X \in ( {\rm dw}\widetilde{\mathcal{A}} )\mbox{}^\perp$ and $A \in {\rm dw}\widetilde{\mathcal{A}}$. Note that $A \in {\rm dw}\widetilde{\mathcal{A}} \cap \mathcal{E}$ since $X \in ( {\rm dw}\widetilde{\mathcal{A}} )\mbox{}^\perp \subseteq \mathcal{E}$, $Y \in \mathcal{E}$ and $\mathcal{E}$ is thick. It follows ${\rm Ext}^1(A,Y) = 0$ and that the previous sequence splits. We have that $Y$ is a direct summand of $X \in ( {\rm dw}\widetilde{\mathcal{A}} )\mbox{}^\perp$, and hence $Y \in ( {\rm dw}\widetilde{\mathcal{A}} )\mbox{}^\perp$ since $( {\rm dw}\widetilde{\mathcal{A}} )\mbox{}^\perp$ is closed under direct summands. 
\end{proof}

The previous result implies that $({\rm dw}\widetilde{\mathcal{P}_n}, ({\rm dw}\widetilde{\mathcal{P}_n})^\perp)$ and $({\rm ex}\widetilde{\mathcal{P}_n}, ({\rm ex}\widetilde{\mathcal{P}_n})^\perp)$ are compatible, and so $({\rm ex}\widetilde{\mathcal{P}_n}, ({\rm ex}\widetilde{\mathcal{P}_n})^\perp)$ is complete, provided $R$ is a Noetherian ring. Using the Hovey's correspondence, Theorem \ref{nprojgrad} follows. Note in its statement that we are not assuming $R$ Noetherian. In fact, it is possible to give a proof of this result that does not depend on this fact. Such a proof consists of constructing for every complex in ${\rm ex}\widetilde{\mathcal{P}_n}$ a filtration by the set $({\rm ex}\widetilde{\mathcal{P}_n})^{\leq \kappa}$, for a fixed infinite cardinal $\kappa > {\rm Card}(R)$. 

Given a projective module $P$ written as a direct sum decomposition of countably generated projective modules $P = \bigoplus_{i \in I} P_i$, note that $P$ is $\kappa$-small if, and only if, ${\rm Card}(I) \leq \kappa$. \\

\begin{definition} \
\begin{itemize}
\item[{\bf (1)}] Given $M \in \mathcal{P}_n$ with a projective resolution \[ {\rm (\ast)} = \left( 0 \longrightarrow \bigoplus_{i \in I_n} P^i_n  \longrightarrow \cdots \longrightarrow \bigoplus_{i \in I_1} P^i_1 \longrightarrow \bigoplus_{i \in I_0} P^i_0 \longrightarrow M \longrightarrow 0 \right). \] We shall say that a projective resolution \[ {\rm (\ast \ast)} = \left( 0 \longrightarrow \bigoplus_{i \in I'_n} P^i_n \longrightarrow \cdots \longrightarrow \bigoplus_{i \in I'_1} P^i_1 \longrightarrow \bigoplus_{i \in I'_0} P^i_0 \longrightarrow N \longrightarrow 0 \right) \] is a \underline{nice sub-resolution} of {\rm ($\ast$)} if {\rm ($\ast \ast$)} is a sub-complex of {\rm ($\ast$)} such that $I'_k \subseteq I_k$ for every $0\leq k \leq n$. 
\end{itemize}
Let $\kappa$ be an infinite cardinal satisfying $\kappa > {\rm Card}(R)$.
\begin{itemize}
\item[{\bf (2)}] We shall denote by $\mathcal{P}_n(\kappa)$ the class of $n$-projective modules $M$ with a projective resolution {\rm ($\ast$)} such that ${\rm Card}(I_k) \leq \kappa$ for every $0 \leq k \leq n$. Note that $\mathcal{P}_n(\kappa) \subseteq (\mathcal{P}_n)^{\leq \kappa}$. 

\item[{\bf (3)}] If we consider the resolutions {\rm ($\ast$)} and {\rm ($\ast \ast$)} above, then we say that {\rm ($\ast \ast$)} is a \underline{nice $\kappa$-small sub-resolution} of {\rm ($\ast$)} if each $I'_k$ is a $\kappa$-small subset of $I_k$. \\
\end{itemize}
\end{definition}

\begin{lemma}\label{lemma2} Let $\kappa$ be an infinite regular cardinal satisfying $\kappa > {\rm Card}(R)$. Let $M \in \mathcal{P}_n$ with a projective resolution given by {\rm ($\ast$)}. For every $\kappa$-small sub-module $N \subseteq M$, there exists a nice $\kappa$-small sub-resolution \newpage \[ {\rm (\ast \ast \ast)} = \left( 0 \longrightarrow \bigoplus_{i \in I'_n} P^i_n \longrightarrow \cdots \longrightarrow \bigoplus_{i \in I'_1} P^i_1 \longrightarrow \bigoplus_{i \in I'_0} P^i_0 \longrightarrow N' \longrightarrow 0 \right) \] of {\rm ($\ast$)} such that $N \subseteq N'$. Moreover, if $N$ has a nice $\kappa$-small sub-resolution {\rm ($\ast \ast$)} of $M$, then {\rm ($\ast \ast \ast$)} can be constructed in such a way that {\rm ($\ast \ast$)} is a sub-resolution of {\rm ($\ast \ast \ast$)}. 
\end{lemma}

\begin{proof} Using the zig-zag procedure in a similar way as in Proposition \ref{Kaplansky}, it is possible to construct $N'$ as described in the statement. Now suppose that $N$ has a nice $\kappa$-small sub-resolution \[ {\rm (\ast \ast)} = \left( 0 \longrightarrow \bigoplus_{i \in I'^N_n} P^i_n \longrightarrow \cdots \longrightarrow \bigoplus_{i \in I^N_1} P^i_1 \longrightarrow \bigoplus_{i \in I^N_0} P^i_0 \longrightarrow N \longrightarrow 0 \right) \] of ($\ast$). Take the quotient of ($\ast$) by this resolution and get \[ 0 \longrightarrow \bigoplus_{i \in I_n - I^N_n} P^i_n \longrightarrow \cdots \longrightarrow \bigoplus_{i \in I_1 - I^N_1} P^i_1 \longrightarrow \bigoplus_{i \in I_0 - I^N_0} P^i_0 \longrightarrow \frac{M}{N} \longrightarrow 0. \] Consider the $\kappa$-small sub-module $\left< z + N \right>$, with $z \not\in N$, and apply the zig-zag procedure to get a projective sub-resolution of the previous one, \[ 0 \longrightarrow \bigoplus_{i \in I'_n - I^N_n} P^i_n \longrightarrow \cdots \longrightarrow \bigoplus_{i \in I'_1 - I^N_1} P^i_1 \longrightarrow \bigoplus_{i \in I'_0 - I^N_0} P^i_0 \longrightarrow \frac{N'}{N} \longrightarrow 0, \] where each set $I'_k - I^N_k$ is a $\kappa$-small set. Finally, note that \[ {\rm (\ast \ast \ast)} = \left( 0 \longrightarrow \bigoplus_{i \in I'_n} P^i_n \longrightarrow \cdots \longrightarrow \bigoplus_{i \in I'_1} P^i_1 \longrightarrow \bigoplus_{i \in I'_0} P^i_0 \longrightarrow N' \longrightarrow 0 \right) \] is a nice $\kappa$-small sub-resolution of ($\ast$) such that it contains $(\ast \ast)$.
\end{proof}

\

\begin{lemma} Let $\kappa$ be an infinite regular cardinal satisfying $\kappa > {\rm Card}(R)$. Let $X \in {\rm dw}\widetilde{\mathcal{P}_n}$ and $Y$ be a $\kappa$-small and bounded above sub-complex of $X$. Then there exists a (bounded above) sub-complex $Y'$ of $X$ such that $Y \subseteq Y'$ and $Y' \in {\rm dw}\widetilde{\mathcal{P}_n(\kappa)}$. 
\end{lemma}
\begin{proof} Since $Y$ is a bounded above complex, there exists $m \in \mathbb{Z}$ such that $Y_k = 0$ for every $k > m$. We are given the following commutative diagram where the vertical arrows are monomorphisms: \newpage
\[ \begin{tikzpicture}
\matrix (m) [matrix of math nodes, row sep=1.5em, column sep=3em]
{ Y =\mbox{ } \cdots & 0 & Y_m & Y_{m-1} & \cdots \\ X =\mbox{ } \cdots & X_{m+1} & X_m & X_{m-1} & \cdots \\ };
\path[->]
(m-1-1) edge (m-1-2)
(m-1-2) edge (m-1-3) edge (m-2-2)
(m-1-3) edge node[above] {$\partial_m$} (m-1-4) edge (m-2-3)
(m-1-4) edge (m-1-5) edge (m-2-4)
(m-2-1) edge (m-2-2) (m-2-2) edge node[above] {$\partial_{m+1}$} (m-2-3) (m-2-3) edge node[above] {$\partial_m$} (m-2-4) (m-2-4) edge (m-2-5);
\end{tikzpicture} \]
Since $X_m$ is an $n$-projective module, we have a projective resolution \[ 0 \longrightarrow \bigoplus_{i \in I_n(m)} P^i_n(m) \longrightarrow \cdots \longrightarrow \bigoplus_{i \in I_1(m)} P^i_1(m) \longrightarrow \bigoplus_{i \in I_0(m)} P^i_0(m) \longrightarrow X_m \longrightarrow 0. \] By the previous lemma, there exists a sub-module $Y'_m$ of $X_m$ containing $Y_m$, along with a nice $\kappa$-small sub-resolution \[ 0 \longrightarrow \bigoplus_{i \in I'_n(m)} P^i_n(m) \longrightarrow \cdots \longrightarrow \bigoplus_{i \in I'_1(m)} P^i_1(m) \longrightarrow \bigoplus_{i \in I'_0(m)} P^i_0(m) \longrightarrow Y'_m \longrightarrow 0. \] Note that ${\rm Card}(\partial_m(Y'_m) + Y_{m-1}) \leq \kappa$ and $Y_{m-1} \subseteq \partial_m(Y'_m) + Y_{m-1} \subseteq X_{m-1}$. Now choose a sub-module $Y'_{m-1} \subseteq X_{m-1}$ such that $\partial_m(Y'_m) + Y_{m-1} \subseteq Y'_{m-1}$ and $Y'_{m-1}$ has a nice $\kappa$-small sub-resolution of a fixed resolution of $X_{m-1}$. Repeat this process infinitely many times in order to obtain a chain complex $Y' = \cdots \longrightarrow 0 \longrightarrow Y'_m \longrightarrow Y'_{m-1} \longrightarrow \cdots$ in ${\rm dw}\widetilde{\mathcal{P}_n(\kappa)}$ such that $Y \subseteq Y' \subseteq X$. 
\end{proof}

\

\begin{theorem} Every complex $X \in {\rm ex}\widetilde{\mathcal{P}_n}$ has a ${\rm ex}\widetilde{\mathcal{P}_n(\kappa)}$-filtration. 
\end{theorem}
\begin{proof} Let $m \in \mathbb{Z}$ be arbitrary and $T_1 \subseteq X_m$ be a $\kappa$-small sub-module of $X_m$. Fix for $X_m$ a projective resolution of length $n$. By Lemma \ref{lemma2} there exists a $\kappa$-small sub-module $Y^1_m$ of $X_m$ such that $T_1 \subseteq Y^1_m$ and that $Y^1_m$ has a nice $\kappa$-small projective sub-resolution of the given resolution of $X_m$. Note that $\partial_m(Y^1_m)$ is a $\kappa$-small sub-module of $X_{m-1}$, so there exists a sub-module $Y^1_{m-1}$ of $X_{m-1}$ such that $\partial_m(Y^1_m) \subseteq Y^1_{m-1}$ and that $Y^1_{m-1}$ has a nice $\kappa$-small projective sub-resolution of the given resolution of $X_{m-1}$. Keep repeating this argument infinitely many times. We obtain a sub-complex $Y^1 = ( \cdots \longrightarrow 0 \longrightarrow Y^1_m \longrightarrow Y^1_{m-1} \longrightarrow \cdots ) \subseteq X$ in ${\rm dw}\widetilde{\mathcal{P}_n(\kappa)}$. Note that $Y^1$ is not necessarily exact. We shall construct a complex $X^1$ from $Y^1$ such that $X^1 \subseteq X$ and $X^1 \in {\rm ex}\widetilde{\mathcal{P}_n(\kappa)}$. The rest of this proof uses an argument similar to the one appearing in \cite[Theorem 4.6]{Rada}. Fix any $p \in \mathbb{Z}$. Then ${\rm Card}(Y^1_p) \leq \kappa$ and so ${\rm Card}(Z_p (Y^1)) \leq \kappa$. Since $X$ is exact and ${\rm Card}(Z_p (Y^1)) \leq \kappa$, there exists a sub-module $U \subseteq X_{p+1}$ with ${\rm Card}(U) \leq \kappa$ such that $Z_p (Y^1) \subseteq \partial_{p+1}(U)$. Let $C^1$ be a $\kappa$-small sub-complex of $X$ such that $U \subseteq C_{p+1}$, $C_j = 0$ for every $j > p+1$, and that each $C_j$ with $j \leq p$ has a nice $\kappa$-small projective sub-resolution of a given resolution of $X_j$. Since $Y^1 + C$ is a bounded above sub-complex of $X$, by the previous lemma there exists a $\kappa$-small sub-complex $Y^2$ of $X$ such that $Y^1 + C \subseteq Y^2$ and that each $Y^2_j$ has a nice $\kappa$-small projective sub-resolution of a given resolution of $X_j$. Note that $Z_p (Y^1) \subseteq \partial_{p+1}(Y^2_{p+1})$. Construct $Y^3$ from $Y^2$ as we just constructed $Y^2$ from $Y^1$, and so on, making sure to use the same $p \in \mathbb{Z}$ at each step. Set $X^1 = \bigcup_{j = 1}^\infty Y^j \subseteq X$. Note that $X^1$ is exact at $p$. Repeat this argument to get exactness at any level. So we may assume that $X^1$ is an exact complex.  Every $X^1_k$ has a nice $\kappa$-small projective sub-resolution of the given resolution of $X_k$: \[ 0 \longrightarrow \bigoplus_{i \in I_{n}(k)} P^i_n(k) \longrightarrow \cdots \longrightarrow \bigoplus_{i \in I_1(k)} P^i_1(k) \longrightarrow \bigoplus_{i \in I_0(k)} P^i_0(k) \longrightarrow X_k \longrightarrow 0. \] Indeed, for every $j$ one has a projective sub-resolution of the form \[ 0 \longrightarrow \bigoplus_{i \in I^j_{n}(k)} P^i_n(k) \longrightarrow \cdots \longrightarrow \bigoplus_{i \in I^j_1(k)} P^i_1(k) \longrightarrow \bigoplus_{i \in I^j_0(k)} P^i_0(k) \longrightarrow Y^j_k \longrightarrow 0, \] where $I^1_l(k) \subseteq I^2_l(k) \subseteq \cdots$ for every $0\leq l \leq n$, by Lemma \ref{lemma2}. If we take the union of all of the previous sequences, then we obtain the following exact sequence: \[ 0 \longrightarrow \bigoplus_{i \in \bigcup_{j\geq i} I^j_{n}(k)} P^i_n(k) \longrightarrow \cdots \longrightarrow \bigoplus_{i \in \bigcup_{j \geq 1} I^j_0(k)} P^i_0(k) \longrightarrow \bigcup_{j \geq 1} Y^j_k = X^1_j \longrightarrow 0, \] where $\bigcup_{j \geq 1} I^j_l(k) \subseteq I_l(k)$ for every $0 \leq l \leq n$. Therefore, $X^1\in {\rm ex}\widetilde{\mathcal{P}_n(\kappa)}$. Moreover, note that the quotient $X / X^1$ belongs to ${\rm ex}\widetilde{\mathcal{P}_n}$. The rest of the proof follows by using transfinite induction. 
\end{proof}

It follows by Proposition \ref{gentransex} that $({\rm ex}\widetilde{\mathcal{P}_n}, ({\rm ex}\widetilde{\mathcal{P}_n})^\perp)$ is a cotorsion pair cogenerated by the set ${\rm ex}\widetilde{\mathcal{P}_n(\kappa)}$, and so it is complete by the Eklof and Trlifaj Theorem. 

Degreewise $n$-projective model structures are not monoidal in general. However, we shall see in the next section that the degreewise projective model structure is monoidal if $R$ is a commutative ring with weak dimension at most $1$. 

We conclude this section presenting the dual of the degreewise $n$-projective model structure. Consider the class $\mathcal{I}_n$ of $n$-injective modules. \\

\begin{theorem}[The \underline{degreewise $n$-injective model structure}] \label{ninjgrad} There exists a unique Abelian model structure on $\Cadl$ where the weak equivalences are the quasi-isomorphisms, the (trivial) fibrations are the epimorphisms with kernels in ${\rm dw}\widetilde{\mathcal{I}_n}$ (in ${\rm ex}\widetilde{\mathcal{I}_n}$), and the (trivial) cofibrations are the monomorphisms with cokernels in $\mbox{}^\perp({\rm ex}\widetilde{\mathcal{I}_n})$ (in $\mbox{}^\perp({\rm dw}\widetilde{\mathcal{I}_n})$). \\
\end{theorem}

Given a Grothendieck category $\mathcal{C}$ with a generator $G$, consider the set $\mathcal{S}$ of all disk complexes $D^m(C)$ with $C \in \Omega^n(G/J)$ and $J$ running over the set of subobjects of $G$. If $D^m(C) \in \mathcal{S}$ and $Y \in {\rm dw}\widetilde{\mathcal{I}_n(\mathcal{C})}$, using \cite[Lemma 3.1 (5)]{Gillespie2} we have ${\rm Ext}^1_{\Complexes}(D^m(C), Y) \cong {\rm Ext}^1_{\mathcal{C}}(C,Y_m)$. On the other hand, ${\rm Ext}^1_{\mathcal{C}}(C,Y_m) \cong {\rm Ext}^{n+1}_{\mathcal{C}}(G/J, Y_m) = 0$, since $Y_m$ is $n$-injective. Then, $D^m(C) \in \mbox{}^\perp( {\rm dw}\widetilde{\mathcal{I}_n(\mathcal{C})} )$. This implies ${\rm dw}\widetilde{\mathcal{I}_n(\mathcal{C})} \subseteq \mathcal{S}^\perp$. The other inclusion follows in the same way, and so ${\rm dw}\widetilde{\mathcal{I}_n(\mathcal{C})} = \mathcal{S}^\perp$. It follows $(\mbox{}^\perp({\rm dw}\widetilde{\mathcal{I}_n}), {\rm dw}\widetilde{\mathcal{I}_n})$ is complete. By Propositions \ref{degreecomp} and \ref{compacogen} and the Hovey's correspondence, the previous theorem follows.


\section{$n$-flat and degreewise $n$-flat model structures}

For $n > 0$, let $\mathcal{F}_n$ denote the class of left $n$-$\mathcal{F}_0$-modules of $\Modl$, where $\mathcal{F}_0$ denotes the class of flat modules. We shall say that $M \in {\rm Ob}(\Modl)$ is \underline{$n$-flat} if $M \in \mathcal{F}_n$. It is known that $(\mathcal{F}_n, (\mathcal{F}_n)^\perp)$ is a hereditary and complete  cotorsion pair (See \cite[Theorem 4.1.3]{Gobel}). It follows by \cite[Corollary 3.8 and Theorem 3.12]{Gillespie2} that we have two compatible cotorsion pairs $(\widetilde{\mathcal{F}_n}, (\widetilde{\mathcal{F}_n})^\perp)$ and $({\rm dg}\widetilde{\mathcal{F}_n}, ({\rm dg}\widetilde{\mathcal{F}_n})^\perp)$. On the other hand, by \cite[Propositions 3.2 and 3.3]{Gillespie}, we have two cotorsion pairs $({\rm dw}\widetilde{\mathcal{F}_n}, ({\rm dw}\widetilde{\mathcal{F}_n})\mbox{}^\perp)$ and $({\rm ex}\widetilde{\mathcal{F}_n}, ({\rm ex}\widetilde{\mathcal{F}_n})\mbox{}^\perp)$. The goal of this section is to construct filtrations for the classes $\widetilde{\mathcal{F}_n}$ and ${\rm ex}\widetilde{\mathcal{F}_n}$ by the sets $(\widetilde{\mathcal{F}_n})^{\leq \kappa}$ and $({\rm ex}\widetilde{\mathcal{F}_n})^{\leq \kappa}$, respectively, based on early investigations developed by S. T. Aldrich, E. Enochs, J. R. Garc\'ia Rozas and L. Oyonarte in \cite{Aldrich1} for the case $n = 0$. After this we deduce the completeness of the previous pairs, and the existence of the following two model structures: \\

\begin{theorem}[The \underline{$n$-flat model structure}] \label{modelonplano} There exists a unique Abelian model structure on $\Cadl$ where the (trivial) cofibrations are the monomorphisms with cokernels in ${\rm dg}\widetilde{\mathcal{F}_n}$ (in $\widetilde{\mathcal{F}_n}$), the (trivial) fibrations are the epimorphisms with kernels in $( \widetilde{\mathcal{F}_n} )\mbox{}^\perp$ (in $( {\rm dg}\widetilde{\mathcal{F}_n} )\mbox{}^\perp$), and the weak equivalences are the quasi-isomorphisms. \\
\end{theorem}

\begin{theorem}[The \underline{degreewise $n$-flat model structure}] \label{degnflatmodel} There exists a unique Abelian model structure on $\Cadl$ where the (trivial) cofibrations are the monomorphisms with cokernels in ${\rm dw}\widetilde{\mathcal{F}_n}$ (in ${\rm ex}\widetilde{\mathcal{F}_n}$), the (trivial) fibrations are  the epimorphisms with kernels in $({\rm ex}\widetilde{\mathcal{F}_n})\mbox{}^\perp$ (in $({\rm dw}\widetilde{\mathcal{F}_n})\mbox{}^\perp$), and the weak equivalences are the quasi-isomorphisms. \\
\end{theorem}

For the case $n = 0$, the completeness of $(\widetilde{\mathcal{F}_0}, (\widetilde{\mathcal{F}_0})^\perp)$ and $({\rm dg}\widetilde{\mathcal{F}_0}, ({\rm dg}\widetilde{\mathcal{F}_0})^\perp)$ was proven by J. Gillespie in \cite{Gillespie2}, using the notions of pure and dg-pure sub-complexes. The compatibility of these two pairs and the Hovey's correspondence allowed Gillespie to find a unique monoidal Abelian model structure on $\Cadl$, called the \underline{flat model structure}, where the weak equivalences are the quasi-isomorphisms, the (trivial) cofibrations are the injections with (exact) dg-flat cokerels, and the (trivial) fibrations are the surjections with (exact) dg-cotorsion kernels (See \cite[Corollary 5.1]{Gillespie2}). \newpage

The cotorsion pair $(\mathcal{F}_0, (\mathcal{F}_0)^\perp)$ was proven to be complete by E. Enochs. His argument consists in showing that every flat module is filtered by the set of $\kappa$-small flat modules, with $\kappa$ an infinite regular cardinal satisfying $\kappa > {\rm Card}(R)$. Enoch's arguments allow us to construct filtrations for the class $\mathcal{F}_n$ by the set $(\mathcal{F}_n)^{\leq \kappa}$. Regarding this matter, the following lemma shall help us to construct the filtrations for $\widetilde{\mathcal{F}_n}$ and ${\rm ex}\widetilde{\mathcal{F}_n}$ mentioned above. \\

\begin{definition} We recall that a sub-module $S$ of a module $M$ is \underline{pure} if the sequence $0 \longrightarrow N \otimes_R S \longrightarrow N \otimes_R M$ is exact, for every right $R$-module $N$. 

If $M$ is an $n$-flat module and $N \subseteq M$ is a sub-module, we say that $N$ is an \underline{$n$-pure sub-module of $M$} if there exist an exact sequence \[ 0 \longrightarrow F_n \longrightarrow \cdots \longrightarrow F_1 \longrightarrow F_0 \longrightarrow M \longrightarrow 0 \] with each $F_k$ flat, and a commutative diagram with exact rows 
\[ \begin{tikzpicture}
\matrix (m) [matrix of math nodes, row sep=1em, column sep=1.5em]
{ 0 & S_n & \cdots & S_1 & S_0 & N & 0 \\ 0 & F_n & \cdots & F_1 & F_0 & M & 0 \\ };
\path[->]
(m-1-1) edge (m-1-2) (m-1-2) edge (m-1-3) edge (m-2-2) (m-1-3) edge (m-1-4) (m-1-4) edge (m-1-5) edge (m-2-4) (m-1-5) edge (m-1-6) edge (m-2-5) (m-1-6) edge (m-1-7) edge (m-2-6)
(m-2-1) edge (m-2-2) (m-2-2) edge (m-2-3) (m-2-3) edge (m-2-4) (m-2-4) edge (m-2-5) (m-2-5) edge (m-2-6) (m-2-6) edge (m-2-7); 
\end{tikzpicture} \]
such that each $S_k \longrightarrow F_k$ is an inclusion with $S_k$ a pure sub-module of $F_k$. \\
\end{definition}

\begin{lemma}\label{puresubres} Let $\kappa$ be an infinite regular cardinal satisfying $\kappa > {\rm Card}(R)$, $M \in \mathcal{F}_n$ with a flat resolution \[ (1) = ( 0 \longrightarrow F_n \stackrel{f_n}\longrightarrow F_{n-1} \longrightarrow \cdots \longrightarrow F_1 \stackrel{f_1}\longrightarrow F_0 \stackrel{f_0}\longrightarrow M \longrightarrow 0 ) \] and $N$ be a $\kappa$-small sub-module of $M$. Then there exists an $n$-pure sub-module $N'$ of $M$ containing $N$ with a resolution \[ 0 \longrightarrow S'_n \longrightarrow \cdots \longrightarrow S'_1 \longrightarrow S'_0 \longrightarrow N' \longrightarrow 0 \] such that $S'_k$ a $\kappa$-small and pure sub-module of $F_k$, for every $0 \leq k \leq n$. Moreover, if $N$ is an $n$-pure sub-module of $M$ with a resolution \[ 0 \longrightarrow S_n \longrightarrow \cdots \longrightarrow S_1 \longrightarrow S_0 \longrightarrow N \longrightarrow 0 \] as in the previous definition, then the above resolution of $N'$ can be constructed in such a way that it contains the given resolution of $N$. 
\end{lemma}
\begin{proof} For every $x \in N$ there exists $y_x \in F_0$ such that $x = f_0(y_x)$. Consider the set $Y := \{ y_x \mbox{ : } x \in N \mbox{ and } f_0(y_x) = x \}$ and the sub-module $\left< Y \right> \subseteq F_0$. Since $\left< Y \right>$ is $\kappa$-small, there exists a $\kappa$-small pure sub-module $S_0(1) \subseteq F_0$ such that $\left< Y \right> \subseteq S_0(1)$, by \cite[Lemma 5.3.12]{Enochs}. Note that $f_0(S_0(1)) \supseteq N$. The rest follows by applying the zig-zag procedure like in the proof of Lemma \ref{lemaoyogen}. In the end we have for each $0 \leq k \leq n$ a collection of modules $(S_k(i) \mbox{ : } i \geq 0)$ such that each $S_k(i)$ is a pure sub-module of $F_k$ and such that \[ {\rm (2)} = (0 \longrightarrow S_n \longrightarrow S_{n-1} \longrightarrow \cdots \longrightarrow S_1 \longrightarrow S_0 \longrightarrow Q \longrightarrow 0), \] is exact, where $S_k = \bigcup_{i \geq 0} S_k(i)$ is a $\kappa$-small and pure sub-module of $F_k$ (and so it is flat) and $Q = {\rm CoKer}(f_1|_{S_1}) \subseteq M$. If we take the quotient of (1) by (2), we get a flat resolution of $M/Q$ of length $n$ (Note that the quotient of a flat module by a pure sub-module is also flat). The rest of the statement follows as in Lemma \ref{lemma2}. 
\end{proof}

There is more to say about the class $\mathcal{F}_n$, and it is the fact that the pair $(\mathcal{F}_n, (\mathcal{F}_n)^\perp)$ is perfect. A cotorsion pair $(\mathcal{A,B})$ in an Abelian category $\mathcal{C}$ is said to be \underline{perfect} if $\mathcal{A}$ is a covering class and $\mathcal{B}$ is an enveloping class, i.e. that every module has an $\mathcal{A}$-cover and a $\mathcal{B}$-envelope. A map $f : A \longrightarrow X$ is said to be an \underline{$\mathcal{A}$-cover} of $X$ if $A \in \mathcal{A}$ and if the following conditions are satisfied: 
\begin{itemize}
\item[{\bf (1)}] If $f : A' \longrightarrow X$ is another map with $A' \in \mathcal{A}$, then there exists a map $h : A' \longrightarrow A$ (not necessarily unique) such that $f' = f \circ h$.

\item[{\bf (2)}] If $A' = A$ in {\bf (1)}, then every map $h : A \longrightarrow A$ satisfying the equality $f' = f \circ h$ is an automorphism of $A$. 
\end{itemize}
The notion of \underline{$\mathcal{B}$-envelope} is dual. In \cite[Corollary 5.32]{Gobel}, the authors prove that if $(\mathcal{A, B})$ is a complete cotorsion pair such that $\mathcal{A}$ is closed under direct limits, then $(\mathcal{A,B})$ is perfect. Since the bifunctors ${\rm Tor}^R_i(-,-)$ preserve filtered colimits, it follows that the class $\mathcal{F}_n$ is closed under direct limits, and hence we can conclude that every left $R$-module has an $n$-flat cover (The case $n = 0$ was known as the {\it Flat Cover Conjecture} until 2000 when it was proven to be true by L. Bican, R. El Bashir and E. Enochs in \cite{Bican}.). We shall see that this result is also valid in the category $\Cadl$ of chain complexes of modules. First, we need to show that $(\widetilde{\mathcal{F}_n},(\widetilde{\mathcal{F}_n})\mbox{}^\perp)$ is complete. In \cite[Definition 4.1.2]{Garcia}, \underline{flat chain complexes} are defined as those complexes $F \in {\rm Ob}(\Cadl)$ such that $F$ is exact and $Z_m(F)$ is a flat module for every $m \in \mathbb{Z}$; or equivalenty, those $F$ such that the functor $- \overline{\otimes} F : \Cadr \longrightarrow {\rm {\bf Ab}}$ is left exact (See \cite[Proposition 5.1.2]{Garcia}). As we did in Proposition \ref{npicomplexes}, we have that a chain complex $X$ is $n$-flat if, and only if, it is exact and $Z_m(X) \in \mathcal{F}_n$, for every $m \in \mathbb{Z}$.

In \cite[Proposition 3.1]{Aldrich1}, the authors prove that each element of a flat chain complex $F$ is contained in a $\kappa$-small flat sub-complex $L \subseteq F$ such that the quotient $F / L$ is also flat. Using the procedure appearing in the given reference, along with Lemma \ref{puresubres}, the following theorem follows. \\

\begin{theorem}\label{teonplanocompleto} Let $\kappa$ be an infinite regular cardinal satisfying $\kappa > {\rm Card}(R)$. For any $n$-flat complex $X \in \widetilde{\mathcal{F}_n}$ and any element $x \in X$ (i.e. $x \in X_k$ for some $k \in \mathbb{Z}$), there exists an exact sub-complex $L \subseteq X$ such that $x \in L$ and such that $Z_m(L)$ is a $k$-small and $n$-pure sub-module of $Z_m(X)$, for every $m \in \mathbb{Z}$. \\
\end{theorem}

Note that if $L$ is the complex provided by the previous result, we have a short exact sequence $0 \longrightarrow L \longrightarrow X \longrightarrow X / L \longrightarrow 0$ in $\Cadl$. Note that $X / L$ is exact, since $L$ and $X$ are and the class of exact complexes is thick. By Lemma \ref{exactofexact}, we have the exact sequence $0 \longrightarrow Z_n(L) \longrightarrow Z_n(X) \longrightarrow Z_n (X / L) \longrightarrow 0$. It follows that $Z_n (X / L) \cong Z_n (X) / Z_n (L)$, which is $n$-flat since $Z_n (X)$ is $n$-flat and $Z_n (L)$ is an $n$-pure sub-module of $Z_n (X)$. Therefore, $X / L \in \widetilde{\mathcal{F}_n}$. It follows that the pair $(\widetilde{\mathcal{F}_n}, (\widetilde{\mathcal{F}_n})\mbox{}^\perp)$ has a cogenerating set, namely the set of $\kappa$-small $n$-flat complexes $(\widetilde{\mathcal{F}_n})^{\leq \kappa}$, by Propositions \ref{extensionespeques} and \ref{gentransex}. So we have a complete cotorsion pair $( \widetilde{\mathcal{F}_n}, (\widetilde{\mathcal{F}_n})\mbox{}^\perp)$. Since $( {\rm dg}\widetilde{\mathcal{F}_n} \cap \mathcal{E}, (\widetilde{\mathcal{F}_n})\mbox{}^\perp)$ and $( {\rm dg}\widetilde{\mathcal{F}_n}, (\widetilde{\mathcal{F}_n})\mbox{}^\perp \cap \mathcal{E})$ are compatible cotorsion pairs and the former is complete, by Proposition \ref{compacogen} the latter is also complete. Therefore, Theorem \ref{modelonplano} follows. 

For $n = 0$, the flat model structure is monoidal, as proved by J. Gillespie in \cite[Corollary 5.1]{Gillespie2}. For $n > 0$, the $n$-flat model structure is not monoidal in general (The same counterexample used for the projective case works.).   

Let $\overline{{\rm Tor}}_i(-,-)$ be the left derived functors of $- \overline{\otimes} -$. For every fixed chain complex $X \in {\rm Ob}(\Cadr)$, the functor $\overline{{\rm Tor}}_1(X,-)$ preserves direct limits (This is a consequence of \cite[Proposition 4.2.1 5]{Garcia}.). It follows the class $\widetilde{\mathcal{F}_n}$ is closed under direct limits, and hence the pair $(\widetilde{\mathcal{F}_n}, (\widetilde{\mathcal{F}_n})\mbox{}^\perp)$ is perfect by  \cite[Corollary 5.32]{Gobel}. Hence, the following result follows. \\

\begin{corollary} Every chain complex in $\Cadl$ has an $n$-flat cover. \\
\end{corollary}

The first approaches to this result are given in \cite[Sections 4.3 and 4.4]{Garcia} in the case $n = 0$, for complexes over a commutative Noetherian ring with finite Krull dimension. 

We finish this section giving a characterization of pure sub-complexes of flat complexes. We recall from \cite[Definition 4.3]{Gillespie2} that a complex $S$ in $\Cadl$ is a \underline{pure sub-complex} of a complex $X$ if the sequence $0 \longrightarrow Y \overline{\otimes} S \longrightarrow Y \overline{\otimes} X$ is exact, for every $Y \in {\rm Ob}(\Cadr)$. \\

\begin{proposition} $S$ is a pure sub-complex of a flat complex $F$ if, and only if, $S$ is exact and $Z_m(S)$ is a pure sub-module of $Z_m(F)$, for every $m \in \mathbb{Z}$.
\end{proposition}

\begin{proof} Suppose $S$ is a pure sub-complex of ${\color{brown!40!black}{F}}$. Then $S$ is flat by \cite[Lemma 4.7]{Gillespie2}, and so exact. It suffices to show $Z_m(S)$ is a pure sub-module of $Z_m({\color{brown!40!black}{F}})$. Let $M$ be a right $R$-module. Consider the sphere complex $S^0(M)$. Since $S$ is a pure sub-complex of $F$, the sequence $0 \longrightarrow S^0(M) \overline{\otimes} S \longrightarrow S^0(M) \overline{\otimes} F$ is exact. At each $m \in \mathbb{Z}$, we have that $(S^0(M) \otimes X)_m = M \otimes_R X_{m}$, for every complex $X \in {\rm Ob}(\Cadl)$. Recall the boundary map $\partial^{S^0(M) \otimes X}_{m+1}$ is given by $y \otimes x \mapsto y \otimes \partial^X_{m+1}(x)$ on generators. It is easy to see that $M \otimes_R B_m(X) = B_m(S^0(M) \otimes X)$. It follows the maps $M \otimes_R B_m(X) \longrightarrow M \otimes_R X_m$ are injective, so $(S^0(M) \overline{\otimes} X)_m = \frac{(S^0(M) \otimes X)_m}{B_m(S^0(M) \otimes X)} = \frac{M \otimes_R X_m}{M \otimes_R B_m(X)} \cong M \otimes_R \frac{X_m}{B_m(X)}$. Since $S$ and $F$ are exact, we get $(S^0(M) \overline{\otimes} S)_m \cong M \otimes_R Z_{m-1}(S)$ and  $(S^0(M) \overline{\otimes} F)_m \cong M \otimes_R Z_{m-1}(F)$. For every $m \in \mathbb{Z}$, $0 \longrightarrow (S^0(M) \overline{\otimes} S)_m \longrightarrow (S^0(M) \overline{\otimes} F)_m$ is exact, and so is $0 \longrightarrow M \otimes_R Z_m(S) \longrightarrow M \otimes_R Z_m(F)$ by the previous isomorphisms. Hence, $Z_m(S)$ is a pure sub-module of $Z_m(F)$. 

Now suppose $S$ is an exact sub-complex of $F$ such that $Z_m(S)$ is a pure sub-module of $Z_m(F)$. Let $A$ be a complex in $\Cadr$. We want to show the sequence $0 \longrightarrow A \overline{\otimes} S \longrightarrow A \overline{\otimes} F$ is exact. Since $Z_{n-k}(S)$ is a pure sub-module of $Z_{n-k}(F)$, we have $0 \longrightarrow A_k \otimes_R Z_{n-k}(S) \longrightarrow A_k \otimes_R Z_{n-k}(F)$ is exact in $\Modl$ for all $k, n \in \mathbb{Z}$. Since $S$ and $F$ are exact complexes, we obtain the following commutative diagram where the top and the bottom rows are exact (Recall Lemma \ref{exactofexact}):
\[ \begin{tikzpicture}
\matrix (m) [matrix of math nodes, row sep=1.5em, column sep=2em]
{ & 0 & 0 \\ 0 & A_k \otimes_R Z_{n-k}(S) & A_k \otimes_R Z_{n-k}(F) \\ 0 & A_k \otimes_R S_{n-k} & A_k \otimes_R F_{n-k}  \\ 0 & A_k \otimes_R Z_{n-k-1}(S) & A_k \otimes_R Z_{n-k-1}(F) \\ & 0 & 0 \\ };
\path[->]
(m-1-2) edge (m-2-2) (m-1-3) edge (m-2-3) 
(m-2-1) edge (m-2-2) (m-2-2) edge (m-2-3) edge (m-3-2) (m-2-3) edge (m-3-3) 
(m-3-1) edge (m-3-2) (m-3-2) edge (m-3-3) edge (m-4-2) (m-3-3) edge (m-4-3) 
(m-4-1) edge (m-4-2) (m-4-2) edge (m-4-3) 
(m-4-2) edge (m-5-2) (m-4-3) edge (m-5-3);
\end{tikzpicture} \]
The columns of this diagram are also exact since the cycles $Z_{n-k}(S)$ and $Z_{n-k}(F)$ are flat modules. Since the class of short exact sequences is closed under extensions, we have that the sequence $0 \longrightarrow A_k \otimes_R S_{n-k} \longrightarrow A_k \otimes_R F_{n-k}$ is exact. It follows $0 \longrightarrow (A \otimes S)_n \longrightarrow (A \otimes F)_n$ is exact for every $n \in \mathbb{Z}$. Then $0 \longrightarrow A \otimes S \longrightarrow A \otimes F$ is an exact sequence of complexes, and so each $0 \longrightarrow B_n(A \otimes S) \longrightarrow B_n(A \otimes F)$ is exact in $\Modl$. Since the class of short exact sequences is thick, we have that each $0 \longrightarrow \frac{(A \otimes S)_n}{B_n(A \otimes S)} \longrightarrow \frac{(A \otimes F)_n}{B_n(A \otimes F)}$ is also exact in $\Modl$. Hence, $0 \longrightarrow A \overline{\otimes} S \longrightarrow A \overline{\otimes} F$ is exact in $\Cadl$.  
\end{proof}

Consider the class $\mathcal{F}_0$ of flat modules. In \cite{Aldrich1}, it is proven that the classes ${\rm dw}\widetilde{\mathcal{F}_0}$ and ${\rm ex}\widetilde{\mathcal{F}_0}$ are filtered by $({\rm dw}\widetilde{\mathcal{F}_0})^{\leq \kappa}$ and $({\rm ex}\widetilde{\mathcal{F}_0})^{\leq \kappa}$, respectively. It follows the cotorsion pairs $({\rm dw}\widetilde{\mathcal{F}_0}, ({\rm dw}\widetilde{\mathcal{F}_0})\mbox{}^\perp)$ and $({\rm ex}\widetilde{\mathcal{F}_0}, ({\rm ex}\widetilde{\mathcal{F}_0})\mbox{}^\perp)$ are complete. Since these two cotorsion pairs turn out to be compatible by Proposition \ref{degreecomp}, we have a unique Abelian model structure on $\Cadl$, which shall be referred as the \underline{degreewise flat model structure}, where the (trivial) cofibrations are the monomorphisms with cokernels in ${\rm dw}\widetilde{\mathcal{F}_0}$ (in ${\rm ex}\widetilde{\mathcal{F}_0}$), the (trivial) fibrations are  the epimorphisms with kernels in $({\rm ex}\widetilde{\mathcal{F}_0})\mbox{}^\perp$ (in $({\rm dw}\widetilde{\mathcal{F}_0})\mbox{}^\perp$), and the weak equivalences are the quasi-isomorphisms. This model structure was discovered by J. Gillespie in \cite{Gillespie}. 

To extend this model structure to $n$-flat modules, we show that the induced cotorsion pairs $({\rm dw}\widetilde{\mathcal{F}_n}, ({\rm dw}\widetilde{\mathcal{F}_n})\mbox{}^\perp)$ and $({\rm ex}\widetilde{\mathcal{F}_n}, ({\rm ex}\widetilde{\mathcal{F}_n})\mbox{}^\perp)$ are complete, applying a modified version of the zig-zag procedure, which also works to show the completeness of the pair $({\rm dw}\widetilde{\mathcal{P}_n}, ({\rm dw}\widetilde{\mathcal{P}_n})\mbox{}^\perp)$, without assuming $R$ Noetherian, as in Section 5. Since $({\rm dw}\widetilde{\mathcal{F}_n}, ({\rm dw}\widetilde{\mathcal{F}_n})\mbox{}^\perp)$ and $({\rm ex}\widetilde{\mathcal{F}_n}, ({\rm ex}\widetilde{\mathcal{F}_n})\mbox{}^\perp)$ shall be also compatible, the following result shall follow.  \\

\begin{theorem}[The \underline{degreewise $n$-flat model structure}] \label{degnflatmodel} There exists a unique Abelian model structure on $\Cadl$ where the (trivial) cofibrations are the monomorphisms with cokernels in ${\rm dw}\widetilde{\mathcal{F}_n}$ (in ${\rm ex}\widetilde{\mathcal{F}_n}$), the (trivial) fibrations are  the epimorphisms with kernels in $({\rm ex}\widetilde{\mathcal{F}_n})\mbox{}^\perp$ (in $({\rm dw}\widetilde{\mathcal{F}_n})\mbox{}^\perp$), and the weak equivalences are the quasi-isomorphisms. \\
\end{theorem}

In the following lines, we shall focus on proving that every chain complex in ${\rm ex}\widetilde{\mathcal{F}_n}$ is filtered by the set $({\rm ex}\widetilde{\mathcal{F}_n})^{\leq \kappa}$ of $\kappa$-small exact degreewise $n$-flat complexes. The proof of the following theorem is based on an argument given in \cite[Proposition 4.1]{Aldrich1} (where the authors prove the case $n = 0$). \\

\begin{theorem}\label{metodo2} Let $X \in {\rm ex}\widetilde{\mathcal{F}_n}$ and $x \in X$ (i.e. $x \in X_m$ for some $m \in \mathbb{Z}$.). Then there exists a sub-complex $Y$ of $X$ with $x \in Y$ and such that $Y_m$ is a $\kappa$-small and $n$-pure sub-module of $X_m$, for every $m \in \mathbb{Z}$. 
\end{theorem}

\begin{proof} Assume without loss of generality that $x \in X_0$. Consider the sub-module $\left< x \right> \subseteq X_0$. Since $X_0 \in \mathcal{F}_n$ and $\left< x \right>$ is $\kappa$-small, we can embed $\left< x \right>$ into a $\kappa$-small and $n$-pure sub-module $Y^1_0 \subseteq X_0$ (Lemma \ref{puresubres}). We can construct a small and exact sub-complex \[ L^1 := ( \cdots \longrightarrow L^1_2 \longrightarrow L^1_1 \longrightarrow Y^1_0 \longrightarrow \partial_0(Y^1_0) \longrightarrow 0 \longrightarrow \cdots ), \] since $X$ is exact. The fact that $\partial_0(Y^1_0)$ is $\kappa$-small implies there exists a $\kappa$-small and $n$-pure sub-module $Y_{-1}^2 \subseteq X_{-1}$. As above, we can construct a $\kappa$-small and exact sub-complex of the form \[ L^2 := ( \cdots \longrightarrow L^2_2 \longrightarrow L^2_1 \longrightarrow L^2_0 \longrightarrow Y^2_{-1} \longrightarrow \partial_{-1}(Y^2_{-1}) \longrightarrow 0 \longrightarrow \cdots ). \] Note that it is possible to construct $L^2$ containing $L^1$. Now embed $L^2_0$ into a $\kappa$-small and $n$-pure sub-module $Y^3_0 \subseteq X_0$. Again, construct a sub-complex \[ L^3 := ( \cdots \longrightarrow L^3_2 \longrightarrow L^3_1 \longrightarrow Y^3_0 \longrightarrow Y^2_{-1} + \partial_0(Y_0^3) \longrightarrow \partial_{-1}(Y_{-1}^2) \longrightarrow 0 \longrightarrow \cdots ) \] containing $L^2$, which is $\kappa$-small and exact. Now let $Y^4_1$ be a $\kappa$-small and $n$-pure sub-module of $X_1$ containing $L^3_1$, and construct an exact and $\kappa$-small complex $L^4$ containing $L^3$ of the form \[ L^4 := ( \cdots \rightarrow L^4_2 \rightarrow Y^4_1 \rightarrow Y^3_0 + \partial_1(Y_1^4) \rightarrow Y^2_{-1} + \partial_0(Y_0^3) \rightarrow \partial_{-1}(Y^2_{-1}) \rightarrow 0 \rightarrow \cdots ) \] Embed $Y^3_0 + \partial_1(Y^4_0)$ into a $\kappa$-small and $n$-pure sub-module $Y^5_0 \subseteq X_0$, and construct an exact and $\kappa$-small sub-complex \[ L^5 := ( \cdots \longrightarrow L^5_2 \longrightarrow L^5_1 \longrightarrow Y^5_0 \longrightarrow Y^2_{-1} + \partial_0(Y^5_0) \longrightarrow \partial_{-1}(Y^2_{-1}) \longrightarrow 0 \longrightarrow \cdots ) \] containing $L^4$. In a similar way, construct $\kappa$-small and exact complexes 
\begin{align*}
L^6 & := ( \cdots \longrightarrow L^6_1 \longrightarrow L^6_0 \longrightarrow Y^6_{-1} \longrightarrow \partial_{-1}(Y^6_{-1}) \longrightarrow 0 \longrightarrow \cdots ) \mbox{ \ and} \\ 
L^7 & := ( \cdots \longrightarrow L^7_1 \longrightarrow L^7_0 \longrightarrow L^7_{-1} \longrightarrow Y^7_{-2} \longrightarrow \partial_{-2}(Y_{-2}^7) \longrightarrow 0 \longrightarrow \cdots )
\end{align*} 
such that $Y_{-1}^6$ is a $\kappa$-small and $n$-pure sub-module of $X_{-1}$ containing $Y^6_{-1} + \partial_0(Y^5_0)$, and $Y_{-2}^7$ is a $\kappa$-small and $n$-pure sub-module of $X_{-2}$ containing $\partial_{-1}(Y_{-1}^6)$. 

Keep repeating this procedure and set $Y := \bigcup_{n \geq 1} L^n$, where $Y_i := \bigcup_{n \geq 1} (L^n)_i$. It is clear that $Y$ is an exact complex. We check that each $Y_m$ is a $\kappa$-small and $n$-pure sub-module of $X_m$. It suffices to show the case $m = 0$. We have \[ Y_0 = Y^1_0 \cup L^2_0 \cup Y^3_0 \cup (Y^3_0 + \partial_1(Y^4_1)) \cup Y^5_0 \cup \cdots = Y^1_0 \cup Y^3_0 \cup Y^5_0 \cup \cdots \] is $\kappa$-small. Consider a left flat resolution of length $n$, \[ (1) = ( 0 \longrightarrow F_n \longrightarrow \cdots \longrightarrow F_1 \longrightarrow F_0 \longrightarrow X_0 \rightarrow 0 ). \] By Lemma \ref{puresubres}, we can construct a sub-resolution \[ 0 \longrightarrow S^1_n \rightarrow \cdots \longrightarrow S^1_1 \longrightarrow S^1_0 \longrightarrow Y^1_0 \longrightarrow 0, \] where $\left< x \right> \subseteq Y^1_0$, and each $S^1_k$ is a $\kappa$-small and pure sub-module of $F_k$. As we did in the previous case, applying Lemma \ref{puresubres} infinitely many times, we can get an ascending chain of sub-resolutions
\begin{align*}
0 & \longrightarrow S^1_n \longrightarrow \cdots \longrightarrow S^1_1 \longrightarrow S^1_0 \longrightarrow Y^1_0 \longrightarrow 0, \\
0 & \longrightarrow S^3_n \longrightarrow \cdots \longrightarrow S^3_1 \longrightarrow S^3_0 \longrightarrow Y^3_0 \longrightarrow 0, \\
0 & \longrightarrow S^5_n \longrightarrow \cdots \longrightarrow S^5_1 \longrightarrow S^5_0 \longrightarrow Y^5_0 \longrightarrow 0, \mbox{ $\dots$}
\end{align*}
Taking the union of these sub-resolutions yields an exact sequence \[ (2) = \left( 0 \longrightarrow \bigcup_j S^j_n \longrightarrow \cdots \longrightarrow \bigcup_j S^j_1 \longrightarrow \bigcup_j S^j_0 \longrightarrow Y_0 \longrightarrow 0 \right), \] where each $\bigcup_j S^j_k$ is a $\kappa$-small and pure sub-module of $F_k$ (so it is flat). 
\end{proof}

\

If $Y$ is the sub-complex of $X$ given by the previous theorem, then we have that $Y_m \in (\mathcal{F}_n)^{\leq \kappa}$ and $X_m / Y_m \in \mathcal{F}_n$ (Take the quotient of (1) by (2)). It follows $Y \in ({\rm ex}\widetilde{\mathcal{F}_n})^{\leq \kappa}$ and $X / Y \in {\rm ex}\widetilde{\mathcal{F}_n}$. \newpage

\begin{remark} The arguments given in the previous proof also work to show that every exact degreewise $n$-projective complex has a filtration by the set ${\rm ex}\widetilde{\mathcal{P}_n(\kappa)}$. Proceed in the same way as in the previous proof, using Lemma \ref{lemma2}, until constructing $Y$. Consider \[ (1) = \left( 0 \longrightarrow \bigoplus_{i \in I_n} P^i_n \longrightarrow \cdots \longrightarrow \bigoplus_{i \in I_1} P^i_1 \longrightarrow \bigoplus_{i \in I_0} P^i_0 \longrightarrow X_0 \longrightarrow 0 \right) \] a projective resolution of $X_0$ of length $n$, where each direct sum consists of countably generated projective modules. By Lemma \ref{lemma2}, we can construct $Y^1_0$ containing $\left< x \right>$ with a sub-resolution of the form \[ (2) = \left( 0 \longrightarrow \bigoplus_{i \in I^1_n} P^i_n \longrightarrow \cdots \longrightarrow \bigoplus_{i \in I^1_1} P^i_1 \longrightarrow \bigoplus_{i \in I^1_0} P^i_0 \longrightarrow Y^1_0 \longrightarrow 0 \right), \] where each $I^1_k \subset I_k$ is a $\kappa$-small. Note that the quotient of {\rm (1)} by {\rm (2)} yields a projective resolution of $X_0 / Y^1_0$ of length $n$, so $X_0 / Y^1_0 \in \mathcal{P}_n$. Using Lemma \ref{lemma2} again, we can construct a sub-resolution containing {\rm (2)}, say \[ (3) = \left( 0 \longrightarrow \bigoplus_{i \in I^3_n} P^i_n \longrightarrow \cdots \longrightarrow \bigoplus_{i \in I^3_1} P^i_1 \longrightarrow \bigoplus_{i \in I^3_0} P^i_0 \longrightarrow Y^3_0 \longrightarrow 0 \right) \] such that $X_0 / Y^3_0 \in \mathcal{P}_n$. We keep applying Lemma \ref{lemma2} to get an ascending chain of sub-resolutions of {\rm (1)}: 
\begin{align*}
0 & \longrightarrow \bigoplus_{i \in I^1_n} P^i_n \longrightarrow \cdots \longrightarrow \bigoplus_{i \in I^1_1} P^i_1 \longrightarrow \bigoplus_{i \in I^1_0} P^i_0 \longrightarrow Y^1_0 \longrightarrow 0, \\
0 & \longrightarrow \bigoplus_{i \in I^3_n} P^i_n \longrightarrow \cdots \longrightarrow \bigoplus_{i \in I^3_1} P^i_1 \longrightarrow \bigoplus_{i \in I^3_0} P^i_0 \longrightarrow Y^3_0 \longrightarrow 0, \\
0 & \longrightarrow \bigoplus_{i \in I^5_n} P^i_n \longrightarrow \cdots \longrightarrow \bigoplus_{i \in I^5_1} P^i_1 \longrightarrow \bigoplus_{i \in I^5_0} P^i_0 \longrightarrow Y^5_0 \longrightarrow 0, \mbox{ \ }\cdots
\end{align*}

Now we take the union of this ascending chain and get the exact complex
\begin{align*}
(4) & = \left( 0 \longrightarrow \bigoplus_{i \in \bigcup_j I^j_n} P^i_n \longrightarrow \cdots \longrightarrow \bigoplus_{i \in \bigcup_j I^j_1} P^i_1 \longrightarrow \bigoplus_{i \in \bigcup_j I^j_0} P^i_0 \longrightarrow Y_0 \longrightarrow 0 \right)
\end{align*}
Since each $\bigcup_j I^j_k$ is a $\kappa$-small subset of $I_k$, we have that the previous sequence is a nice $\kappa$-small projective sub-resolution of {\rm (1)}. Note also that the quotient of {\rm (1)} by {\rm (4)} yields a projective resolution of $X_0 / Y_0$ of length $n$. Then $Y_0 \in \mathcal{P}_n(\kappa)$. In a similar way, we have that $Y_m \in \mathcal{P}_n(\kappa)$ and $X_m / Y_m \in \mathcal{P}_n$, for every $m \in \mathbb{Z}$. 
\end{remark}

By the previous result we have that $({\rm ex}\widetilde{\mathcal{P}_n}, ({\rm ex}\widetilde{\mathcal{P}_n})\mbox{}^\perp)$ and $({\rm ex}\widetilde{\mathcal{F}_n}, ({\rm ex}\widetilde{\mathcal{F}_n})\mbox{}^\perp)$ are complete cotorsion pairs. Using a similar argument one can prove the same for the pairs $({\rm dw}\widetilde{\mathcal{P}_n}, ({\rm dw}\widetilde{\mathcal{P}_n})\mbox{}^\perp)$ and $({\rm dw}\widetilde{\mathcal{F}_n}, ({\rm dw}\widetilde{\mathcal{F}_n})\mbox{}^\perp)$. By Proposition \ref{degreecomp}, the pairs $({\rm dw}\widetilde{\mathcal{F}_n}, ({\rm dw}\widetilde{\mathcal{F}_n})\mbox{}^\perp)$ and $({\rm ex}\widetilde{\mathcal{F}_n}, ({\rm ex}\widetilde{\mathcal{F}_n})\mbox{}^\perp)$ are compatible. Therefore, Theorem \ref{degnflatmodel} follows. Similarly, we have gotten another proof of Theorem \ref{nprojgrad}. 

J. Gillespie showed in \cite[Subsection 5.2]{Gillespie} that the degreewise flat model structure is not monoidal on $(\Cadl, \otimes)$ in general. He considered the category of complexes over the ring $\mathbb{Z}_4$, where $Y = \cdots \longrightarrow \mathbb{Z}_4 \stackrel{\times 2}\longrightarrow \mathbb{Z}_4 \stackrel{\times 2}\longrightarrow \mathbb{Z}_4 \longrightarrow \cdots$ is an exact degreewise flat complex, but $Y \otimes Y \not\in {\rm ex}\widetilde{\mathcal{F}_0}$, since it is not even exact. This counterexample also works to show that the degreewise $n$-flat model structure is not monoidal for every $n > 0$. 

One problem with the ring $\mathbb{Z}_4$ considered by Gillespie is that $\mathbb{Z}_4$ has infinite weak dimension. For $\mathbb{Z}_2$ is a $\mathbb{Z}_4$-module with infinite flat dimension, since ${\rm Tor}^{\mathbb{Z}_4}_k(\mathbb{Z}_2, \mathbb{Z}_2) \cong \mathbb{Z}_2 \neq 0$ for every $k \geq 0$ \cite[Chapter 3, Example 9]{Osborne}. This "inconvenient" can be solved if we impose an extra condition on $R$. \\ 

\begin{proposition}\label{degflatmono} Let $R$ be a commutative ring with weak dimension at most $1$. Then the degreewise projective and degreewise flat model structures on $\Cadl$ are monoidal. \\
\end{proposition}

We know that $(\Cadl, \otimes)$ is a closed symmetric monoidal category with unit $S^0(R)$. Flat objects in a monoidal category $(\mathcal{C},\otimes)$ are defined as those $X \in {\rm Ob}(\mathcal{C})$ such that the functor $- \otimes X: \mathcal{C} \longrightarrow \mathcal{C}$ is exact. Previously, we have considered flat objects in $\Cadl$ with respect to $\overline{\otimes}$. For the rest of this section, we shall say that a complex $F$ is \underline{flat with respect to $\otimes$} if the the functor $- \otimes F$ is exact. 

We only prove the previous proposition for the flat case (since the projective case follows similarly). By \cite[Theorem 4.2]{Hovey2007} it suffices to verify the following conditions:
\begin{itemize}
\item[{\bf (1)}] Every complex in ${\rm dw}\widetilde{\mathcal{F}_0}$ is flat with respect to $\otimes$.
\item[{\bf (2)}] If $X, Y \in {\rm dw}\widetilde{\mathcal{F}_0}$, then so is $X \otimes Y$.
\item[{\bf (3)}] If $X, Y \in {\rm dw}\widetilde{\mathcal{F}_0}$ and one of them is exact, then $X \otimes Y$ is exact.
\item[{\bf (4)}] $S^0(R) \in {\rm dw}\widetilde{\mathcal{F}_0}$. 
\end{itemize}

Note that condition {\bf (4)} is immediate, and {\bf (2)} is easy to verify. To show {\bf (1)}, we establish the following characterization. \\ 

\begin{proposition} A chain complex $X \in \Cadl$ over a commutative ring $R$ is flat with respect to $\otimes$ if, and only if, it is degreewise flat. 
\end{proposition}
\begin{proof} Let $Y$ be a chain complex such that $- \otimes Y$ is exact. Consider an exact sequence $0 \longrightarrow A \longrightarrow B \longrightarrow C \longrightarrow 0$ in $\Modl$. Then we have an exact sequence $0 \longrightarrow S^0(A) \longrightarrow S^0(B) \longrightarrow S^0(C) \longrightarrow 0$ in $\Cadl$. It follows $0 \longrightarrow S^0(A) \otimes Y \longrightarrow S^0(B) \otimes Y \longrightarrow S^0(C) \otimes Y \longrightarrow 0$ is also exact. So for each $n \in \mathbb{Z}$ we have the exact sequence $0 \longrightarrow A \otimes_R Y_n \longrightarrow B \otimes_R Y_n \longrightarrow C \otimes_R Y_n \longrightarrow 0$. 

Now suppose $Y \in {\rm dw}\widetilde{\mathcal{F}_0}$. Consider a short exact sequence of chain complexes $0 \longrightarrow A \stackrel{\alpha}\longrightarrow B \stackrel{\beta}\longrightarrow C \longrightarrow 0$ and apply $- \otimes Y$. We need to check that for every $n \in \mathbb{Z}$, the sequence $0 \longrightarrow (A \otimes Y)_n \stackrel{(\alpha \otimes Y)_n}\longrightarrow (B \otimes Y)_n \stackrel{(\beta \otimes Y)_n}\longrightarrow (C \otimes Y)_n \longrightarrow 0$ is exact. In other words, we shall see \\
\[ \begin{tikzpicture}[descr/.style={fill=white}]
\matrix (m) [matrix of math nodes, row sep=1em, column sep=6em, text height=1.5ex, text depth=0.25ex]
{ (\ast) \mbox{ \ = \ } (\mbox{ } 0 & \bigoplus_{k \in \mathbb{Z}} A_k \otimes_R Y_{n-k} & \bigoplus_{k \in \mathbb{Z}} B_k \otimes_R Y_{n-k} \\ & \bigoplus_{k \in \mathbb{Z}} C_k \otimes_R Y_{n-k} & 0 \mbox{ }) \\ };
\path[overlay,->, font=\scriptsize, >=latex]
(m-1-1) edge (m-1-2) (m-1-2) edge node[above] {$\bigoplus_{k \in \mathbb{Z}} \alpha_k \otimes_R Y_{n-k}$} (m-1-3) 
(m-1-3) edge [out=355, in=175] node[descr] {$\bigoplus_{k \in \mathbb{Z}} \beta_k \otimes_R Y_{n-k}$} (m-2-2)
(m-2-2) edge (m-2-3);
\end{tikzpicture} \] 
is exact. For every $k \in \mathbb{Z}$, the sequence \[ 0 \longrightarrow A_k \otimes_R Y_{n-k} \stackrel{\alpha_k \otimes_R Y_{n-k}}\longrightarrow B_k \otimes_R Y_{n-k} \stackrel{\beta_k \otimes Y_{n-k}}\longrightarrow C_k \otimes_R Y_{n-k} \longrightarrow 0 \] is exact since $Y_{n-k}$ is flat. It follows ($\ast$) is exact since the direct sum of exact sequences is exact (Recall that homology commutes with direct sums.). 
\end{proof}

\

It is only left to verify {\bf (3)}. If $X$ and $Y$ are degreewise flat chain complexes over a commutative ring $R$ with weak dimension at most $1$, then by \cite[Theorem 9.16]{Osborne} we have the following K\"unneth exact sequence for every $k \in \mathbb{Z}$: 
\[ \begin{tikzpicture}[descr/.style={fill=white}]
\matrix (m) [matrix of math nodes, row sep=1em, column sep=2em, text height=1.5ex, text depth=0.25ex]
{ 0 & \bigoplus_{i+j = k} H_i(X) \otimes H_j(Y) & H_k(X \otimes Y) \\ & \bigoplus_{i+j = k-1} {\rm Tor}^R_1(H_i(X), H_j(Y)) & 0 \\ };
\path[overlay,->, font=\scriptsize, >=latex]
(m-1-1) edge (m-1-2) (m-1-2) edge (m-1-3) 
(m-1-3) edge [out=355, in=175] (m-2-2)
(m-2-2) edge (m-2-3);
\end{tikzpicture} \]   
If in addition $X$ or $Y$ is exact, then it follows that $X \otimes Y$ is also exact. Therefore, Proposition \ref{degflatmono} follows.


\section*{Acknowledgements} 

The author would like to thank Professor Andr\'e Joyal of the Department of Mathematics at UQ\`AM, for his support and suggestions, and his enlightening introduction to the theory of modules over ringoids. This work has been financially supported by scholarships from ISM-CIRGET and Fondation de l'UQ\`AM. 


\end{document}